%% file: Hilbert-polygons.tex
\documentclass[12pt]{amsart}
\usepackage{amsmath} 
\usepackage{amsfonts} 
\usepackage{amssymb} 
\usepackage{amsthm} 
\usepackage{amscd} 
\usepackage{enumerate}
\usepackage{graphicx}
\numberwithin{equation}{section}

\input{./latex-macros}

\newcommand{\bC}{\partial \cC} 
\newcommand{\bP}{\partial \cP} 
\newcommand{\bQ}{\partial \cQ} 
\newcommand{\bS}{\partial \cS} 
\newcommand{\bT}{\partial \cT} 
\newcommand{\dC}{d_{\cC}} 
\newcommand{\dP}{d_{\cP}}

\newcommand{\FC}{F_{\cC}} 
\newcommand{\FP}{F_{\cP}}
\newcommand{\FQ}{F_{\cQ}}
\newcommand{\FS}{F_{\cS}}
\newcommand{\FT}{F_{\cT}}
\renewcommand{\t}{\tau}
\renewcommand{\inc}{\subseteq} 
\renewcommand{\diam}[2]{\mathrm{diam}_{#1} \! \left( #2 \right)} 

\begin{document}

\title[]{Hilbert geometry for convex polygonal domains}

\author{Bruno Colbois}
\address{Bruno Colbois, 
Universit\'{e} de Neuch\^{a}tel, 
Institut de math\'{e}matique, 
Rue \'{E}mile Argand~11, 
Case postale~158, 
CH--2009 Neuch\^{a}tel, 
Switzerland}
\email{bruno.colbois@unine.ch}

\author{Constantin Vernicos}
\address{Constantin Vernicos, 
UMR 5149 du CNRS \& Universit\'{e} Montpellier~I\kern -0.05em I, 
Institut de math\'{e}\-matique et de mod\'{e}lisation de Montpellier, 
Place Eug\`{e}ne Bataillon, 
Case courrier~51, \linebreak 
F--34095 Montpellier Cedex, 
France}
\email{constantin.vernicos@math.univ-montp2.fr}

\author{Patrick Verovic}
\address{Patrick Verovic, 
UMR 5127 du CNRS \& Universit\'{e} de Savoie, 
Laboratoire de math\'{e}matique, 
Campus scientifique, 
F--73376 Le Bourget-du-Lac Cedex, 
France}
\email{verovic@univ-savoie.fr}

\date{\today}
\subjclass[2000]{Primary: global Finsler geometry, Secondary: convexity}


\begin{abstract}
We prove in this paper that the Hilbert geometry associated with an open convex polygonal 
set is Lipschitz equivalent to Euclidean plane. 
\end{abstract}

\maketitle

\bigskip
\bigskip


\section{Introduction} 

A \emph{Hilbert domain} in $\Rn{m}$ is a metric space $(\cC , \dC)$, where $\cC$ is an 
\emph{open bounded convex} set in $\Rn{m}$ and $\dC$ is the distance function on $\cC$ 
--- called the \emph{Hilbert metric} --- defined as follows. 
   
\medskip
   
Given two distinct points $p$ and $q$ in $\cC$, 
let $a$ and $b$ be the intersection points of the straight line defined by $p$ and $q$ 
with $\bC$ so that $p = (1 - s) a + s b$ and $q = (1 - t) a + t b$ with $0 < s < t < 1$. 
Then 
$$
\dC(p , q) \as \frac{1}{2} \ln{\! [a , p , q , b]},
$$ 
where 
$$
[a , p , q , b] \as \frac{1 - s}{s} \times \frac{t}{1 - t} > 1
$$ 
is the cross ratio of the $4$-tuple of ordered collinear points $(a , p , q , b)$ 
(see Figure~\ref{fig:Hilbert-metric}). 

We complete the definition by setting $\dC(p , p) \as 0$. 

\begin{figure}[h]
   \includegraphics[width=9.7cm,height=6cm,keepaspectratio=true]{./Figures/polygons-fig-1.eps}
   \caption{\label{fig:Hilbert-metric} The Hilbert metric $\dC$}
\end{figure}

\bigskip

The metric space $(\cC , \dC)$ thus obtained is a complete non-compact geodesic metric space 
whose topology is the one induced by the canonical topology of $\Rn{m}$ and in which the affine open segments 
joining two points of the boundary $\bC$ are geodesics that are isometric to $(\RR , |\cdot|)$. 
It is to be mentioned here that in general the affine segment between two points in $\cC$ 
may \emph{not} be the \emph{unique} geodesic joining these points (for example, if $\cC$ is a square).  Nevertheless, this uniqueness holds whenever $\cC$ is \emph{strictly} convex. 

\bigskip

For further information about Hilbert geometry, we refer to \cite{Ben06,Bus55,BusKel53,Egl97,Gol88,Ver05} 
and the excellent introduction \cite{Soc00} by Soci\'{e}-M\'{e}thou. 

\bigskip
\bigskip

The two fundamental examples of Hilbert domains $(\cC , \dC)$ in $\Rn{m}$ correspond to the case 
when\,\,{}$\cC$ is an ellipsoid, which gives the Klein model of $m$-dimensional hyperbolic geometry 
(see for example \cite[first chapter]{Soc00}), and the case when\,\,{}$\clos{\cC}$ is a $m$-simplex, 
for which there exists a norm $\norm{\cdot}_{\cC}$ on $\Rn{m}$ such that $(\cC , \dC)$ 
is isometric to the normed vector space $(\Rn{m} , \norm{\cdot}_{\cC})$ 
(see \cite[pages 110--113]{dlH93} or \cite[pages 22--23]{Nus88}). 
Therefore, it is natural to study the Hilbert domains $(\cC , \dC)$ in $\Rn{m}$ for which $\cC$ 
is close to either an ellipsoid or a $m$-simplex. 

\medskip

The first and last authors thus proved in \cite{ColVero04} that any Hilbert domain $(\cC , \dC)$ in $\Rn{m}$ 
such that the boundary $\bC$ is a $\Cl{2}$ hypersurface with non-vanishing Gaussian curvature 
is Lipschitz equivalent to $m$-dimensional hyperbolic space $\Hn{m}$. 

\medskip

On the other hand, F\"{o}rtsch and Karlsson showed in \cite{ForKar05} that a Hilbert domain in $\Rn{m}$ 
is isometric to a normed vector space if and only if it is given by a $m$-simplex. 
In addition, Lins established in his PhD thesis \cite[Lemma~2.2.5]{Lin07} that 
the Hilbert geometry associated with an open convex polygonal set in $\Rn{2}$ 
can be isometrically embedded in the normed vector space $(\Rn{N^{2}} , \supnorm{\cdot})$, 
where $N$ is the number of vertices of the polygon. 

\bigskip
\bigskip

The aim of this paper is to prove that the Hilbert geometry associated with 
an open convex \emph{polygonal} set $\cP$ in $\Rn{2}$ 
is Lipschitz equivalent to Euclidean plane (Theorem~\ref{thm:bi-Lipschitz} in the last section). 
A straighforward consequence of this result is that \emph{all} the Hilbert \emph{polygonal} domains 
in $\Rn{2}$ are Lipschitz equivalent to each other, 
which is a fact that is far from being obvious at a first glance. 

\smallskip

The idea of the proof is to decompose a given open convex $n$-sided polygon $\cP$ 
into $n$ triangles having one common vertex in $\cP$ and whose opposite edges to that vertex 
are the sides of $\cP$, and then to show that each of these triangles is Lipschitz equivalent to the cone 
it defines with that vertex. This second point is the most technical part of the paper 
and is based on Proposition~\ref{prop:comparison-1}. 

\bigskip

\begin{remark*} 
It seems that our result might be extended 
to higher dimensions to prove more generally that any Hilbert domain in $\Rn{m}$ given by a \emph{polytope} 
is Lipschitz equivalent to $m$-dimensional Euclidean space. 
Nevertheless, computations in that case appear to be much more difficult 
since they involve not only the edges of the polytope but also its faces. 
\end{remark*}

\bigskip
\bigskip
\bigskip


\section{Preliminaries} \label{sec:preliminaries} 

This section is devoted to some technical properties we will need 
for the proof of Theorem~\ref{thm:bi-Lipschitz} in Section~\ref{sec:Lipschitz-eq}. 
The key results are contained in Proposition~\ref{prop:map} and Proposition~\ref{prop:comparison-1}. 

\bigskip

Let us first recall that the distance function $\dC$ is associated with the Finsler metric 
$\FC$ on $\cC$ given, for any $p \in \cC$ and any $v \in T_{\! p}\cC = \Rn{m}$ 
(tangent vector space to $\cC$ at $p$), by 
$$
\FC(p , v) 
\as 
\frac{1}{2} \! \left( \frac{1}{t^{-}} + \frac{1}{t^{+}} \right) \ \ \mbox{if} \ \ v \neq 0,
$$ 
where $t^{-} = t_{\cC}^{-}(p , v)$ and $t^{+} = t_{\cC}^{+}(p , v)$ are the \emph{unique positive} numbers 
such that $p - t^{-} v \in \bC$ and $p + t^{+} v \in \bC$, 
and $\FC(p , 0) \as 0$ (see Figure~\ref{fig:Finsler-metric}). 

\medskip

This means that for every $p , q \in \cC$ and $v \in T_{\! p}\cC = \Rn{m}$, we have 
$\disp \FC(p , v) = \der{}{t}_{{^{\big |}}_{t = 0}} \!\!\!\!\!\! \dC(p , p + t v)$ 
and $\dC(p , q)$ is the infimum of the length $\disp \int{0}{1}{\FC(\s(t) , \s'(t))}{t}$ 
with respect to $\FC$ when $\s : [0 , 1] \to \cC$ ranges over all the $\Cl{1}$ curves joining $p$ to $q$. 

\medskip

\begin{remark*} 
For $p \in \cC$ and $v \in T_{\! p}\cC = \Rn{m}$ with $v \neq 0$, we will define 
$p^{-} = p_{\cC}^{-}(p , v) \as p - t_{\cC}^{-}(p , v) v$ 
and $p^{+} = p_{\cC}^{+}(p , v) \as p + t_{\cC}^{+}(p , v) v$. 
Then, given any arbitrary norm $\norm{\cdot}$ on $\Rn{m}$, we can write 
$$
\FC(p , v) 
= 
\frac{1}{2} \norm{v} \!\! \left( \frac{1}{\norm{p - p^{-}}} + \frac{1}{\norm{p - p^{+}}} \right) \! .
$$ 
\end{remark*}

\begin{figure}[h]
   \includegraphics[width=9.7cm,height=6cm,keepaspectratio=true]{./Figures/polygons-fig-2.eps}
   \caption{\label{fig:Finsler-metric} The Finsler metric $\FC$}
\end{figure}

\bigskip

\begin{figure}[h]
   \includegraphics[width=10cm,height=10cm,keepaspectratio=true]{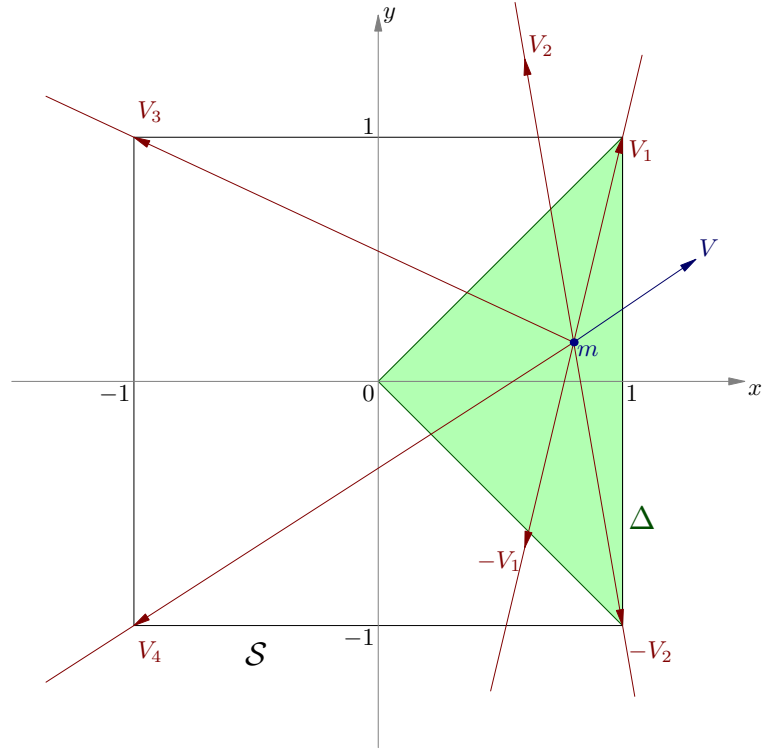}
   \caption{\label{fig:zones} The six zones for the vector $V$}
\end{figure}

\begin{notations}

Let $\cS \as ]{-1} , 1[ \: \cart \: ]{-1} , 1[ \; \inc \Rn{2}$ be the standard open square, 
$\D \as \{ (x , y) \in \Rn{2} \st |y| < x < 1 \} \inc \cS$ 
the open triangle whose vertices are $0 = (0 , 0)$, $(1 , -1)$ and $(1 , 1)$, 
and $\cZ \as \{ (X , Y) \in \Rn{2} \st |Y| < X \} \inc \Rn{2}$ the open cone associated with $\D$ 
(see Figure~\ref{fig:zones}). 

\bigskip

The canonical basis of $\Rn{2}$ will be denoted by $(e_{1} , e_{2})$. 

\bigskip

The usual $\ll{1}$-norm on $\Rn{2}$ and its associated distance will be denoted respectively 
by $\norm{\cdot}$ and $d$. 

\end{notations}

\bigskip

\begin{definition} \label{def:sector} 
   For any pair $(V_{1} , V_{2})$ of vectors in $\Rn{2} \setmin \{ 0 \}$, the set 
   $$
   S(V_{1} , V_{2}) \as \{ s V_{1} + t V_{2} \st s \geq 0 \ \ \mbox{and} \ \ t \geq 0 \}
   $$ 
   will be called the \emph{sector} associated with this pair. 
\end{definition}

\bigskip

\begin{remark*}
The sector $S(V_{1} , V_{2})$ is the convex hull of the set $(\RR_{+} V_{1}) \cup (\RR_{+} V_{2})$. 
\end{remark*}

\bigskip

Let us begin with the following useful lemma: 

\begin{lemma} \label{lem:sector} 
   Given a basis $(V_{1} , V_{2})$ of $\Rn{2}$ and a vector $V \in \Rn{2}$, we have 
   $$
   V \in S(V_{1} , V_{2}) 
   \iff 
   \left( \detb{(V_{1} , V_{2})}{(V_{1} , V)} \geq 0 
   \quad \mbox{and} \quad 
   \detb{(V_{1} , V_{2})}{(V , V_{2})} \geq 0 \right).
   $$ 
\end{lemma}

\medskip

\begin{proof}~\\ 
The lemma is a mere consequence of the fact that the coordinate system $(s , t)$ 
of any vector $V$ in $\Rn{2}$ with respect to a basis $(V_{1} , V_{2})$ of $\Rn{2}$ is equal to 
$\disp \left( \detb{(V_{1} , V_{2})}{(V , V_{2})} \ , \ \detb{(V_{1} , V_{2})}{(V_{1} , V)} \right)$. 
\end{proof}

\bigskip

Now, we have 

\begin{proposition} \label{prop:map} 
   The map $\F : \cS \to \Rn{2}$ defined by 
   $$
   \F(x , y) = (X , Y) \as (\atanh(x) \ , \ \atanh(y))
   $$ 
   is a smooth diffeomorphism such that 
   
   \begin{enumerate}
      \item $\F(\D) = \cZ$, and 
      
      \smallskip
      
      \item for all $m \in \D$ and $V \in T_{\! m}\cS = \Rn{2}$, 
      ~$\FS(m , V) \leq \norm{\Tg{m}{\F}{V}} \leq 2 \FS(m , V)$. 
   \end{enumerate}
\end{proposition}

\bigskip

Before proving this result, we will need the following (see Figure~\ref{fig:zones}): 

\begin{lemma} \label{lem:zones} 
   Let $m = (x , y) \in \D \inc \cS$, and define in $T_{\! m}\cS = \Rn{2}$ the vectors 
   
   \smallskip
   
   \begin{center} 
   $V_{1} \as (1 , 1) - m = (1 - x \ , \ 1 - y)$, 
   \quad 
   $V_{2} \as m - (1 , -1) = (-1 + x \ , \ 1 + y)$, 
   
   \medskip
   
   $V_{3} \as (-1 , 1) - m = (-1 - x \ , \ 1 - y)$ 
   \quad and \quad 
   $V_{4} \as (-1 , -1) - m = (-1 - x \ , \ -1 - y)$. 
   \end{center}
   
   \medskip
   
   Then we have the inclusions 
   
   \medskip
   
   \begin{enumerate}
      \item $\disp S(V_{1} , V_{2}) \inc \left\{ V = (\l , \m) \in \Rn{2} ~\bigg{|}~ 
      \m > 0 \ \ \mbox{and} \ \ \frac{|\l|}{1 - x^{2}} \leq \frac{\m}{1 - y^{2}} \right\}$, 
      
      \medskip
      
      \item $\disp S(V_{2} , V_{3}) \inc \{ V = (\l , \m) \in \Rn{2} \st \l < 0 < \m \}$, 
      
      \medskip
      
      \item $\disp S(V_{3} , V_{4}) \inc \left\{ V = (\l , \m) \in \Rn{2} ~\bigg{|}~ 
      \l < 0 \ \ \mbox{and} \ \ \frac{|\m|}{1 - y^{2}} \leq \frac{-\l}{1 - x^{2}} \right\}$, and 
      
      \medskip
      
      \item $S(V_{4} , -V_{1}) \inc \{ V = (\l , \m) \in \Rn{2} \st \l < 0 \ \ \mbox{and} \ \ \m < 0 \}$. 
   \end{enumerate}
\end{lemma}

\bigskip

\begin{proof}~\\ 
First of all, we have $\detb{(e_{1} , e_{2})}{(V_{1} , V_{2})} = 2 (1 - x) > 0$, 
$\detb{(e_{1} , e_{2})}{(V_{2} , V_{3})} = 2 (x + y) > 0$, 
$\detb{(e_{1} , e_{2})}{(V_{3} , V_{4})} = 2 (1 - x) > 0$ 
and $\detb{(e_{1} , e_{2})}{(V_{4} , -V_{1})} = 2 (x - y) > 0$ (since $x > y$). 

\smallskip

This shows that $(V_{1} , V_{2})$, $(V_{2} , V_{3})$, $(V_{3} , V_{4})$ and $(V_{4} , -V_{1})$ are 
all bases of $\RR^{2}$ having the same orientation as $(e_{1} , e_{2})$. 

\medskip

Then, let $V = (\l , \m)$ be an arbitrary vector in $T_{\! m}\cS = \Rn{2}$. 

\bigskip

$\bullet$ \textsf{Point~(1):} 
If $V \in S(V_{1} , V_{2})$, then, according to Lemma~\ref{lem:sector}, we have 
$$
0 \leq \detb{(e_{1} , e_{2})}{(V_{1} , V)} = (1 - x) \m - (1 - y) \l 
\quad \mbox{and} \quad 
0 \leq \detb{(e_{1} , e_{2})}{(V , V_{2})} = (1 + y) \l + (1 - x) \m
$$ 
since $(V_{1} , V_{2})$ is a basis of $\RR^{2}$ having the same orientation as $(e_{1} , e_{2})$. 

\smallskip

This writes 
$$
\frac{\l}{1 - x} \leq \frac{\m}{1 - y} 
\quad \mbox{and} \quad 
\frac{-\m}{1 + y} \leq \frac{\l}{1 - x}~,
$$ 
and hence, multiplying both inequalities by $\disp \frac{1}{1 + x} > 0$, we get 
\begin{equation} \label{equ:1-1} 
   \frac{\l}{1 - x^{2}} \leq \frac{\m}{(1 + x) (1 - y)} 
   \quad \mbox{and} \quad 
   \frac{-\m}{(1 + x) (1 + y)} \leq \frac{\l}{1 - x^{2}}~. 
\end{equation}

\smallskip

On the other hand, writing $V = s V_{1} + t V_{2}$ with 
$$
s \as \detb{(V_{1} , V_{2})}{(V , V_{2})} \geq 0 
\quad \mbox{and} \quad 
t \as \detb{(V_{1} , V_{2})}{(V_{1} , V)} \geq 0,
$$ 
the second coordinate $\m$ of $V$ with respect to 
the canonical basis $(e_{1} , e_{2})$ of $\Rn{2}$ equals 
$$
\m = \detb{(e_{1} , e_{2})}{(e_{1} , V)} 
= 
s \detb{(e_{1} , e_{2})}{(e_{1} , V_{1})} + t \detb{(e_{1} , e_{2})}{(e_{1} , V_{2})} 
= s (1 - y) + t (1 + y) > 0.
$$ 

\smallskip

This yields $\disp \frac{\m}{(1 + x) (1 - y)} \leq \frac{\m}{1 - y^{2}}$ 
since $0 < 1 + y \leq 1 + x$, and hence 
\begin{equation} \label{equ:1-2} 
   \frac{\l}{1 - x^{2}} \leq \frac{\m}{1 - y^{2}} 
\end{equation} 
from the first part of Equation~\ref{equ:1-1}. 

\smallskip

Moreover, we also have $\disp \frac{-\m}{1 - y^{2}} \leq \frac{-\m}{(1 + x) (1 + y)}$ 
since $0 < 1 - y \leq 1 + x$. Thus, 
\begin{equation} \label{equ:1-3} 
   \frac{-\m}{1 - y^{2}} \leq \frac{\l}{1 - x^{2}} 
\end{equation} 
from the second part of Equation~\ref{equ:1-1}. 

\smallskip

Finally, summarizing Equations~\ref{equ:1-2} and~\ref{equ:1-3}, we obtain 
$\disp \frac{|\l|}{1 - x^{2}} \leq \frac{\m}{1 - y^{2}}$~. 

\bigskip

$\bullet$ \textsf{Point~(2):} 
If $V \in S(V_{2} , V_{3})$, let us write $V = s V_{2} + t V_{3}$ with 
$$
s \as \detb{(V_{2} , V_{3})}{(V , V_{3})} \geq 0 
\quad \mbox{and} \quad 
t \as \detb{(V_{2} , V_{3})}{(V_{2} , V)} \geq 0.
$$ 

\smallskip

Then the first coordinate $\l$ of $V$ with respect to 
the canonical basis $(e_{1} , e_{2})$ of $\Rn{2}$ equals 
$$
\l = \detb{(e_{1} , e_{2})}{(V , e_{2})} 
= 
s \detb{(e_{1} , e_{2})}{(V_{2} , e_{2})} + t \detb{(e_{1} , e_{2})}{(V_{3} , e_{2})} 
= -s (1 - x) - t (1 + x) < 0.
$$ 

\smallskip

On the other hand, the second coordinate $\m$ of $V$ with respect to $(e_{1} , e_{2})$ is equal to 
$$
\m = \detb{(e_{1} , e_{2})}{(e_{1} , V)} 
= 
s \detb{(e_{1} , e_{2})}{(e_{1} , V_{2})} + t \detb{(e_{1} , e_{2})}{(e_{1} , V_{3})} 
= s (1 + y) + t (1 - y) > 0.
$$ 

\bigskip

$\bullet$ \textsf{Point~(3):} 
If $V \in S(V_{3} , V_{4})$, then, according to Lemma~\ref{lem:sector}, we have 
$$
0 \leq \detb{(e_{1} , e_{2})}{(V_{3} , V)} = -(1 + x) \m - (1 - y) \l 
\quad \mbox{and} \quad 
0 \leq \detb{(e_{1} , e_{2})}{(V , V_{4})} = -(1 + y) \l + (1 + x) \m
$$ 
since $(V_{3} , V_{4})$ is a basis of $\RR^{2}$ having the same orientation as $(e_{1} , e_{2})$. 

\smallskip

This writes 
$$
\frac{\m}{1 - y} \leq \frac{-\l}{1 + x} 
\quad \mbox{and} \quad 
\frac{\l}{1 + x} \leq \frac{\m}{1 + y}~,
$$ 
and hence, multiplying the first inequality by $\disp \frac{1}{1 + y} > 0$ 
and the second one by $\disp \frac{1}{1 - y} > 0$, we get 
\begin{equation} \label{equ:3-1} 
   \frac{\m}{1 - y^{2}} \leq \frac{-\l}{(1 + x) (1 + y)} 
   \quad \mbox{and} \quad 
   \frac{\l}{(1 + x) (1 - y)} \leq \frac{\m}{1 - y^{2}}~. 
\end{equation} 

\smallskip

On the other hand, writing $V = s V_{3} + t V_{4}$ with 
$$
s \as \detb{(V_{3} , V_{4})}{(V , V_{4})} \geq 0 
\quad \mbox{and} \quad 
t \as \detb{(V_{3} , V_{4})}{(V_{3} , V)} \geq 0,
$$ 
the first coordinate $\m$ of $V$ with respect to 
the canonical basis $(e_{1} , e_{2})$ of $\Rn{2}$ equals 
$$
\l = \detb{(e_{1} , e_{2})}{(V , e_{2})} 
= 
s \detb{(e_{1} , e_{2})}{(V_{3} , e_{2})} + t \detb{(e_{1} , e_{2})}{(V_{4} , e_{2})} 
= -s (1 + x) - t (1 + x) < 0.
$$ 

\smallskip

This yields $\disp \frac{-\l}{(1 + x) (1 + y)} \leq \frac{-\l}{1 - x^{2}}$ 
since $0 < 1 - x \leq 1 + y$, and hence 
\begin{equation} \label{equ:3-2} 
   \frac{\m}{1 - y^{2}} \leq \frac{-\l}{1 - x^{2}} 
\end{equation} 
from the first part of Equation~\ref{equ:3-1}. 

\smallskip

Moreover, we also have $\disp \frac{\l}{1 - x^{2}} \leq \frac{\l}{(1 + x) (1 - y)}$ 
since $0 < 1 - x \leq 1 - y$. Thus, 
\begin{equation} \label{equ:3-3} 
   \frac{\l}{1 - x^{2}} \leq \frac{\m}{1 - y^{2}} 
\end{equation} 
from the second part of Equation~\ref{equ:3-1}. 

\smallskip

Finally, summarizing Equations~\ref{equ:3-2} and~\ref{equ:3-3}, we obtain 
$\disp \frac{|\m|}{1 - y^{2}} \leq \frac{-\l}{1 - x^{2}}$~.

\bigskip

$\bullet$ \textsf{Point~(4):} 
If $V \in S(V_{4} , -V_{1})$, let us write $V = s V_{4} - t V_{1}$ with 
$$
s \as \detb{(V_{4} , -V_{1})}{(V , V_{4})} \geq 0 
\quad \mbox{and} \quad 
t \as \detb{(V_{4} , -V_{1})}{(-V_{1} , V)} \geq 0.
$$ 

\smallskip

Then the first coordinate $\l$ of $V$ with respect to 
the canonical basis $(e_{1} , e_{2})$ of $\Rn{2}$ equals 
$$
\l = \detb{(e_{1} , e_{2})}{(V , e_{2})} 
= 
s \detb{(e_{1} , e_{2})}{(V_{4} , e_{2})} - t \detb{(e_{1} , e_{2})}{(V_{1} , e_{2})} 
= -s (1 + x) - t (1 - x) < 0.
$$ 

\smallskip

On the other hand, the second coordinate $\m$ of $V$ with respect to $(e_{1} , e_{2})$ is equal to 
$$
\m = \detb{(e_{1} , e_{2})}{(e_{1} , V)} 
= 
s \detb{(e_{1} , e_{2})}{(e_{1} , V_{4})} - t \detb{(e_{1} , e_{2})}{(e_{1} , V_{1})} 
= -s (1 + y) - t (1 - y) < 0.
$$ 
\end{proof}

\bigskip

\begin{proof}[Proof of Proposition~\ref{prop:map}]~\\ 
Only the second point has to be proved since the first one is obvious. 

\smallskip

So, fix $m = (x , y) \in \D \inc \cS$ and $V = (\l , \m) \in T_{\! m}\cS = \Rn{2}$ such that $V \neq 0$. 

\smallskip

A straightforward computation shows that 
$$
\norm{\Tg{m}{\F}{V}} = \left( \frac{\l}{1 - x^{2}} \ , \ \frac{\m}{1 - y^{2}} \right) \! ,
$$ 
and thus 
$$
\norm{\Tg{m}{\F}{V}} = \frac{|\l|}{1 - x^{2}} + \frac{|\m|}{1 - y^{2}}~.
$$ 

\smallskip

Now, let us define the vectors $V_{1}$, $V_{2}$, $V_{3}$ and $V_{4}$ in $T_{\! m}\cS = \Rn{2}$ 
as in Lemma~\ref{lem:zones}. 

\smallskip

Since $\Rn{2}$ is equal to the union of the sectors $S(V_{1} , V_{2})$, $S(V_{2} , V_{3})$, 
$S(V_{3} , V_{4})$, $S(V_{4} , -V_{1})$ and their images by the symmetry 
about the origin $0$, and since the Finsler metric $\FS$ on $\cS$ is reversible, there are four cases 
to be considered. 

\bigskip

$\bullet$ \textsf{Case~1:} $V \in S(V_{1} , V_{2})$. 

\smallskip

The unique positive numbers $\t^{-}$ and $\t^{+}$ such that 
$m - \t^{-} V \in \bS$ and $m + \t^{+} V \in \bS$ satisfy 
$y - \t^{-} \m = -1$ and $y + \t^{+} \m = 1$. 
So, $\t^{-} = (1 + y) / \m$ and $\t^{+} = (1 - y) / \m$, and hence 
$$
\FS(m , V) = \frac{1}{2} \! \left( \frac{1}{\t^{-}} + \frac{1}{\t^{+}} \right) = \frac{\m}{1 - y^{2}}~.
$$ 

\smallskip

But 
$$
\norm{\Tg{m}{\F}{V}} 
\geq 
\frac{|\m|}{1 - y^{2}} 
= 
\frac{\m}{1 - y^{2}}
$$ 
since $\m > 0$ by point~(1) in Lemma~\ref{lem:zones}. 

\smallskip

Therefore, we have 
$$
\FS(m , V) \leq \norm{\Tg{m}{\F}{V}}.
$$ 

\smallskip

On the other hand, point~(1) in Lemma~\ref{lem:zones} yields 
$$
\norm{\Tg{m}{\F}{V}} 
= 
\frac{|\l|}{1 - x^{2}} + \frac{\m}{1 - y^{2}} 
\leq \frac{2 \m}{1 - y^{2}}~,
$$ 
which shows that 
$$
\norm{\Tg{m}{\F}{V}} \leq 2 \FS(m , V).
$$ 

\bigskip

$\bullet$ \textsf{Case~2:} $V \in S(V_{2} , V_{3})$. 

\smallskip

The unique positive numbers $\t^{-}$ and $\t^{+}$ such that 
$m - \t^{-} V \in \bS$ and $m + \t^{+} V \in \bS$ satisfy 
$x - \t^{-} \l = 1$ and $y + \t^{+} \m = 1$. 
So, $\t^{-} = -(1 - x) / \l$ and $\t^{+} = (1 - y) / \m$, and hence 
$$
\FS(m , V) = \frac{1}{2} \! \left( \frac{-\l}{1 - x} + \frac{\m}{1 - y} \right).
$$ 

\smallskip

But point~(2) in Lemma~\ref{lem:zones} implies 
$$
\norm{\Tg{m}{\F}{V}} 
= 
\frac{-\l}{1 - x^{2}} + \frac{\m}{1 - y^{2}} 
= 
\left( \frac{1}{1 + x} \right) \!\! \left( \frac{-\l}{1 - x} \right) 
+ \left( \frac{1}{1 + y} \right) \!\! \left( \frac{\m}{1 - y} \right)
$$ 
with $\disp \frac{-\l}{1 - x} > 0$ and $\disp \frac{\m}{1 - y} > 0$. 

\smallskip

Therefore, since $\disp \frac{1}{2} \leq \frac{1}{1 + x} \leq 1$ 
and $\disp \frac{1}{2} \leq \frac{1}{1 + y} \leq 1$, we have 
$$
\FS(m , V) \leq \norm{\Tg{m}{\F}{V}} \leq 2 \FS(m , V).
$$ 

\bigskip

$\bullet$ \textsf{Case~3:} $V \in S(V_{3} , V_{4})$. 

\smallskip

The unique positive numbers $\t^{-}$ and $\t^{+}$ such that 
$m - \t^{-} V \in \bS$ and $m + \t^{+} V \in \bS$ satisfy 
$x - \t^{-} \l = 1$ and $x + \t^{+} \l = -1$. 
So, $\t^{-} = -(1 - x) / \l$ and $\t^{+} = -(1 + x) / \l$, and hence 
$$
\FS(m , V) = \frac{-\l}{1 - x^{2}}~.
$$ 

\smallskip

But 
$$
\norm{\Tg{m}{\F}{V}} 
\geq 
\frac{|\l|}{1 - x^{2}} 
= 
\frac{-\l}{1 - x^{2}}
$$ 
since $\l < 0$ by point~(3) in Lemma~\ref{lem:zones}. 

\smallskip

Therefore, we have 
$$
\FS(m , V) \leq \norm{\Tg{m}{\F}{V}}.
$$ 

\smallskip

On the other hand, point~(3) in Lemma~\ref{lem:zones} yields 
$$
\norm{\Tg{m}{\F}{V}} 
= 
\frac{-\l}{1 - x^{2}} + \frac{|\m|}{1 - y^{2}} 
\leq \frac{-2 \l}{1 - x^{2}}~,
$$ 
which shows that 
$$
\norm{\Tg{m}{\F}{V}} \leq 2 \FS(m , V).
$$ 

\bigskip

$\bullet$ \textsf{Case~4:} $V \in S(V_{4} , -V_{1})$. 

\smallskip

The unique positive numbers $\t^{-}$ and $\t^{+}$ such that 
$m - \t^{-} V \in \bS$ and $m + \t^{+} V \in \bS$ satisfy 
$x - \t^{-} \l = 1$ and $y + \t^{+} \m = -1$. 
So, $\t^{-} = -(1 - x) / \l$ and $\t^{+} = -(1 + y) / \m$, and hence 
$$
\FS(m , V) = \frac{1}{2} \! \left( \frac{-\l}{1 - x} + \frac{-\m}{1 + y} \right).
$$ 

\smallskip

But point~(2) in Lemma~\ref{lem:zones} implies 
$$
\norm{\Tg{m}{\F}{V}} 
= 
\frac{-\l}{1 - x^{2}} + \frac{-\m}{1 - y^{2}} 
= 
\left( \frac{1}{1 + x} \right) \!\! \left( \frac{-\l}{1 - x} \right) 
+ \left( \frac{1}{1 + y} \right) \!\! \left( \frac{-\m}{1 - y} \right)
$$ 
with $\disp \frac{-\l}{1 - x} > 0$ and $\disp \frac{-\m}{1 - y} > 0$. 

\smallskip

Therefore, since $\disp \frac{1}{2} \leq \frac{1}{1 + x} \leq 1$ 
and $\disp \frac{1}{2} \leq \frac{1}{1 + y} \leq 1$, we have 
$$
\FS(m , V) \leq \norm{\Tg{m}{\F}{V}} \leq 2 \FS(m , V).
$$ 
\end{proof}

\bigskip

\begin{remark*} 
It is to be pointed out that the Lipschitz constants $1$ and $2$ obtained in Proposition~\ref{prop:map} 
are optimal. Indeed, taking $m \as (1 / 2 \, , \, 0) \in \D$ 
and $V \as (0 , 1) \in S(V_{1} , V_{2}) \inc T_{\! m}\cS = \Rn{2}$, 
we get $\norm{\Tg{m}{\F}{V}} = \FS(m , V)$. On the other hand, we have 
$\disp \frac{\norm{\Tg{m}{\F}{V}}}{\FS(m , V)} \to 2$ when $m \to (0 , 0)$ and $V \to (1 , 1)$ 
with $m \in \D$ and $V \in S(V_{1} , V_{2}) \inc T_{\! m}\cS = \Rn{2}$. 
\end{remark*}

\bigskip

Given real numbers $a \in (0 , 1)$ and $c > b \geq 1$, 
let $\cT \inc \Rn{2}$ be the triangle defined as the open convex hull 
of the points $(1 , -1)$, $(1 , 1)$ and $(-a , 0)$, 
and let $\cQ \inc \Rn{2}$ be the quadrilateral defined as the open convex hull 
of the points $(1 , -1)$, $(1 , 1)$, $(-b , c)$ and $(-b , -c)$ (see Figure~\ref{fig:comparison-1}). 

\begin{figure}[h]
   \includegraphics[width=10cm,height=10cm,keepaspectratio=true]{./Figures/polygons-fig-4.eps}
   \caption{\label{fig:comparison-1} The triangle $\cT$ and the quadrilateral $\cQ$}
\end{figure}

\bigskip

Then we have $\D \inc \cT \inc \cQ$ and 

\begin{proposition} \label{prop:comparison-1} 
   There exists a constant $A = A(a , b , c) \in (0 , 1]$ such that 
   $$
   A \FT(m , V) \leq \FQ(m , V) \leq \FT(m , V)
   $$ 
   for all $m \in \D$ and $V \in T_{\! m}\cT = T_{\! m}\cQ = \Rn{2}$. 
\end{proposition}

\bigskip

\begin{remark*}
This is the key result of this section, but also the most technical one of the paper. 
So, the reader may skip it in a first reading without any loss of keeping track of the ideas 
that lead to the final theorem in Section~\ref{sec:Lipschitz-eq}. 
\end{remark*}

\bigskip

\begin{figure}[h]
   \includegraphics[width=9.7cm,height=6cm,keepaspectratio=true]{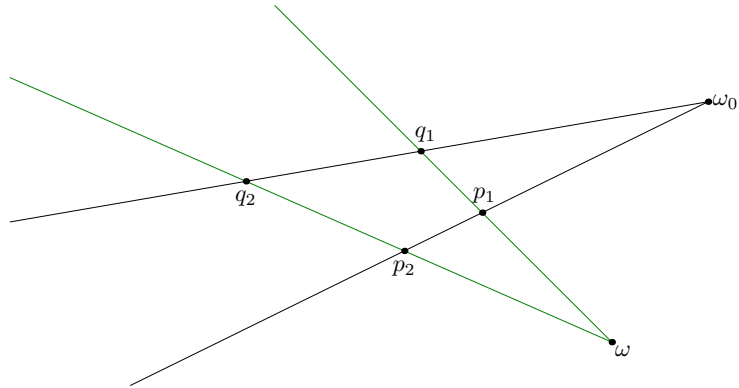}
   \caption{\label{fig:projection} Intersecting pairs of half-lines}
\end{figure}

Before proving Proposition~\ref{prop:comparison-1}, 
we will need the following simple but very useful fact (see Figure~\ref{fig:projection}): 

\begin{lemma} \label{lem:projection} 
   Let $\om$, $q_{1}$ and $q_{2}$ be non-collinear points in $\Rn{2}$, 
   and consider $p_{1} \in \; ]\om , q_{1}[$ and $p_{2} \in \; ]\om , q_{2}[$ such that 
   the lines $(p_{1} p_{2})$ and $(q_{1} q_{2})$ intersect in a point $\om_{0}$. 
   
   If $q_{1} \in \; ]\om_{0} , q_{2}[$, then we have 
   $\disp \frac{\bar{\om q_{2}}}{\bar{\om p_{2}}} > \frac{\bar{\om q_{1}}}{\bar{\om p_{1}}}$~. 
\end{lemma}

\bigskip

\begin{proof}~\\ 
Let $\pi : \Rn{2} \to \Rn{2}$ be the projection of $\Rn{2}$ onto the line $(\om q_{2})$ 
along the direction of $(\om_{0} p_{2})$. 

\medskip

Since $\pi$ is affine, it is barycentre-preserving, and hence $q_{1} \in \; ]\om_{0} , q_{2}[$ 
necessarily implies $\pi(q_{1}) \in \; ]\pi(\om_{0}) \ , \ \pi(q_{2})[ \ = \ ]p_{2} , q_{2}[$. 
So, $\bar{\om q_{2}} > \bar{\om \pi(q_{1})}$. 

\medskip

But, $\pi$ being affine with $\pi(\om) = \om$ and $\pi(p_{1}) = p_{2}$, we also have 
$\disp \frac{\bar{\om \pi(q_{1})}}{\bar{\om p_{2}}} = \frac{\bar{\om q_{1}}}{\bar{\om p_{1}}}$, 
which proves Lemma~\ref{lem:projection} thanks to the previous inequality. 
\end{proof}

\bigskip

\begin{proof}[Proof of Proposition~\ref{prop:comparison-1}]~\\ 
Since $\cT \inc \cQ$, we already have the second inequality. 
So, the very thing to be proved here is the first inequality. 

\medskip

Recall that $\norm{\cdot}$ and $d$ denote respectively the usual $\ll{1}$-norm on $\Rn{2}$ 
and its associated distance. 

\smallskip

Define $\k_{0} \as \diam{d}{\cQ} > 0$, the diameter of $\cQ$ with respect to $d$, 
and let $\theta_{0}$ be the intersection point of the line $\RR (1 , 1)$ 
with the line passing through the points $(1 , -1)$ and $(-b , c)$; 
in other words, $\theta_{0} = (\a_{0} , \a_{0})$ with $\disp \a_{0} \as \frac{c - b}{c + b + 2} \in (0 , 1)$. 

\medskip

Next consider 
$$
\D^{\! +} \as \{ (x , y) \in \Rn{2} \st 0 \leq y < x < 1 \} \inc \D 
\quad \mbox{and} \quad 
\S \as \{ (x , y) \in \Rn{2} \st \a_{0} \leq y < x < 1 \} \inc \D^{\! +},
$$ 
fix $m = (x , y) \in \S$, and define in $T_{\! m}\cT = T_{\! m}\cQ = \Rn{2}$ the vectors 

\smallskip

\begin{center} 
$V_{1} \as (1 , 1) - m = (1 - x \ , \ 1 - y)$, 
\quad 
$V_{2} \as m - (1 , -1) = (-1 + x \ , \ 1 + y)$, 

\medskip

$V_{3} \as (-b , c) - m = (-b - x \ , \ c - y)$ 
\quad and \quad 
$V_{4} \as (-a , 0) - m = (-a - x \ , \ -y)$. 
\end{center}

\smallskip

Then we have 
\begin{eqnarray*} 
   \detb{(e_{1} , e_{2})}{(V_{1} , V_{2})} 
   & = & 
   \detb{(e_{1} , e_{2})}{(-V_{2} \, , V_{1})} 
   = \detb{(e_{1} , e_{2})}{(-V_{1} \, , -V_{2})} = 2 (1 - x) > 0, \\ 
   \detb{(e_{1} , e_{2})}{(V_{2} , V_{3})} 
   & = & 
   \detb{(e_{1} , e_{2})}{(-V_{2} \, , -V_{3})} \\ 
   & = & 
   (1 + c) x + (1 + b) y + b - c \geq (2 + b + c) y + b - c \geq 0 
   \quad \mbox{(since $x \geq y \geq \a_{0}$)}, \\ 
   \detb{(e_{1} , e_{2})}{(V_{1} , V_{3})} 
   & = & 
   \detb{(e_{1} , e_{2})}{(-V_{3} \, , V_{1})} \\ 
   & = & 
   (1 - c) x - (1 + b) y + b + c > (b + c) (1 - x) > 0 
   \quad \mbox{(since $x > y$)}, \\ 
   \detb{(e_{1} , e_{2})}{(V_{3} , V_{4})} 
   & = & 
   (b - a) y + c (x + a) > 0 
   \quad \mbox{(since $b \geq 1 \geq a$)}, \\ 
   \detb{(e_{1} , e_{2})}{(V_{2} , V_{4})} 
   & = & 
   \detb{(e_{1} , e_{2})}{(V_{4} , -V_{2})} 
   = \detb{(e_{1} , e_{2})}{(-V_{2} \, , -V_{4})} = (1 + a) y + x + a > 0, 
   \quad \mbox{and} \\ 
   \detb{(e_{1} , e_{2})}{(V_{1} , V_{4})} 
   & = & 
   \detb{(e_{1} , e_{2})}{(V_{4} , -V_{1})} 
   = \detb{(e_{1} , e_{2})}{(-V_{4} \, , V_{1})} \\ 
   & = & 
   a (1 - y) + x - y > a (1 - y) > 0 
   \quad \mbox{(since $x > y$)}. 
\end{eqnarray*} 

\smallskip

This shows that $(V_{1} , V_{3})$, $(V_{1} , V_{4})$, $(V_{3} , V_{4})$, $(V_{4} , -V_{2})$ 
and $(-V_{2} \, , V_{1})$ are bases of $\RR^{2}$ having the same orientation as $(e_{1} , e_{2})$ 
with $V_{2} \in S(V_{1} , V_{3}) \cap S(V_{1} , V_{4})$, $V_{3} \in S(V_{1} , V_{4})$, 
$-V_{1} \in S(V_{4} , -V_{2})$ and $-V_{3} , -V_{4} \in S(-V_{2} \, , V_{1})$. 

\medskip
   
Given an arbitrary vector $V = (\l , \m) \in T_{\! m}\cT = T_{\! m}\cQ = \Rn{2}$ such that $V \neq 0$, 
there are now four cases to be dealt with. 

\bigskip 

$\bullet$ \textsf{Case~1:} $V \in S(V_{1} , V_{2})$ (see Figure~\ref{fig:case-1}). 

\begin{figure}[h]
   \includegraphics[width=10cm,height=13cm,keepaspectratio=true]{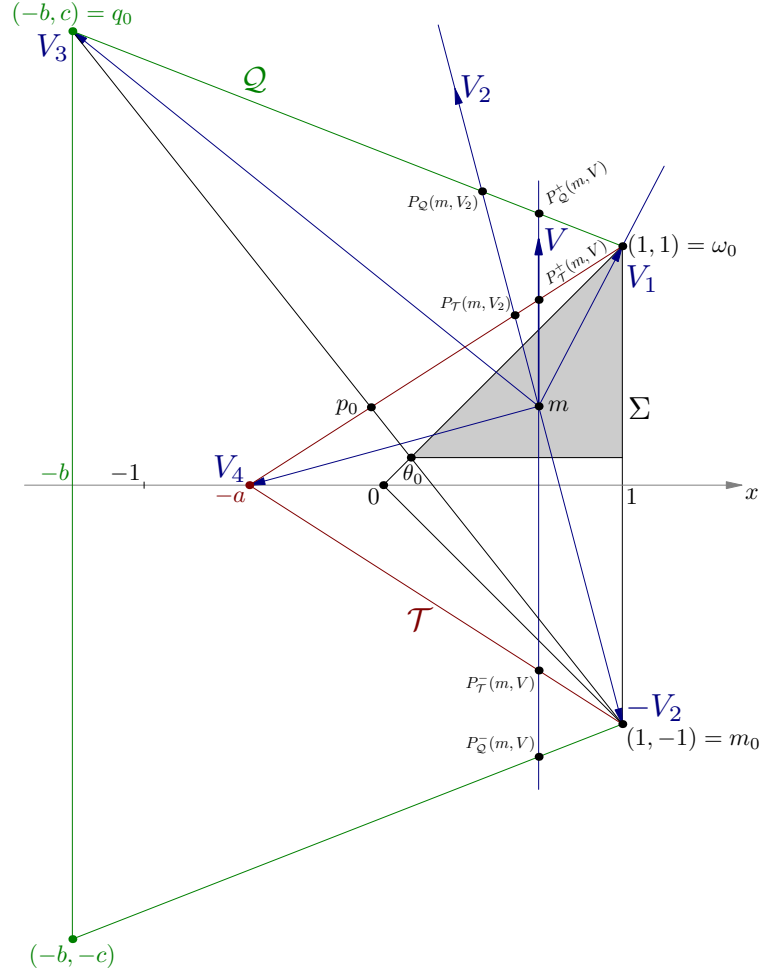}
   \caption{\label{fig:case-1} The case when $V \in S(V_{1} , V_{2})$}
\end{figure}

\bigskip

From $V_{2} \in S(V_{1} , V_{3})$, the half-line $m + \RR_{+} V_{2}$ intersects 
with the segment $[(1 , 1) \, , \, (-b , c)] \inc \bQ$, and hence the same holds for the half-line 
$m + \RR_{+} V$ since $V \in S(V_{1} , V_{2}) \inc S(V_{1} , V_{3})$. 

\smallskip

Moreover, since $V_{2} \in S(V_{1} , V_{4})$, we have $V \in S(V_{1} , V_{2}) \inc S(V_{1} , V_{4})$, 
and this implies that the half-lines $m + \RR_{+} V_{2}$ and $m + \RR_{+} V$ also intersect 
with the segment $[(1 , 1) \, , \, (-a , 0)] \inc \bT$. 

\smallskip

Therefore, if $V \not \in \{ V_{1} , V_{2} \}$, Lemma~\ref{lem:projection} 
with $\om \as m$, $p_{1} \as p_{\cT}^{+}(m , V)$, $q_{1} \as p_{\cQ}^{+}(m , V)$, 
$p_{2} \as p_{\cT}^{+}(m , V_{2})$, $q_{2} \as p_{\cQ}^{+}(m , V_{2})$ and $\om_{0} \as (1 , 1)$ gives 
\begin{equation} \label{equ:inequality-1} 
   \frac{t_{\cQ}^{+}(m , V_{2})}{t_{\cT}^{+}(m , V_{2})} 
   \geq 
   \frac{t_{\cQ}^{+}(m , V)}{t_{\cT}^{+}(m , V)}~, 
\end{equation} 
and this still holds when $V \as V_{1}$ or $V \as V_{2}$. 

\smallskip

On the other hand, if $m_{0} \as (1 , -1)$, $q_{0} \as (-b , c)$ 
and $p_{0}$ is the intersection point of the line passing through $(1 , 1)$ and $(-a , 0)$ 
with the line $(m_{0} q_{0})$, Lemma~\ref{lem:projection} with $\om \as m_{0}$, 
$p_{1} \as p_{\cT}^{+}(m , V_{2})$, $q_{1} \as p_{\cQ}^{+}(m , V_{2})$, 
$p_{2} \as p_{0}$, $q_{2} \as q_{0}$ and $\om_{0} \as (1 , 1)$ yields 
\begin{equation} \label{equ:inequality-2} 
   \frac{\bar{m_{0} q_{0}}}{\bar{m_{0} p_{0}}} 
   \geq 
   \frac{\bar{m_{0} p_{\cQ}^{+}(m , V_{2})}}{\bar{m_{0} p_{\cT}^{+}(m , V_{2})}}~. 
\end{equation} 

\smallskip

But $t_{\cQ}^{+}(m , V_{2}) = \bar{m p_{\cQ}^{+}(m , V_{2})} 
\geq \bar{m p_{\cT}^{+}(m , V_{2})} = t_{\cT}^{+}(m , V_{2}) > 0$ (since $\cT \inc \cQ$) 
and $\bar{m_{0} m} > 0$, which implies 
\begin{equation} \label{equ:inequality-3} 
   \frac{\bar{m_{0} p_{\cQ}^{+}(m , V_{2})}}{\bar{m_{0} p_{\cT}^{+}(m , V_{2})}} 
   \geq 1 \geq 
   \frac{\bar{m_{0} p_{\cT}^{+}(m , V_{2})}}{\bar{m_{0} p_{\cQ}^{+}(m , V_{2})}} 
   = 
   \frac{\bar{m_{0} m} + t_{\cT}^{+}(m , V_{2})}{\bar{m_{0} m} + t_{\cQ}^{+}(m , V_{2})} 
   \geq 
   \frac{t_{\cT}^{+}(m , V_{2})}{t_{\cQ}^{+}(m , V_{2})}~. 
\end{equation} 

\smallskip

Then, combining Equations~\ref{equ:inequality-1}, \ref{equ:inequality-2} and~\ref{equ:inequality-3}, we get 
\begin{equation} \label{equ:c1-t+} 
   \frac{\bar{m_{0} q_{0}}}{\bar{m_{0} p_{0}}} 
   \geq 
   \frac{t_{\cQ}^{+}(m , V)}{t_{\cT}^{+}(m , V)}~. 
\end{equation} 

\medskip

Furthermore, since $m \in \S$ 
and $p_{\cT}^{-}(m , V) \in [(-a , 0) \, , \, (1 , -1)] \inc \RR \cart (-\infty , 0]$ 
(indeed, $-V_{1} \in S(V_{4} , -V_{2})$ implies $-V \in S(-V_{1} \, , -V_{2}) \inc S(V_{4} , -V_{2})$), 
we have $\norm{m - p_{\cT}^{-}(m , V)} \geq d(m , \RR \cart \{ 0 \}) \geq \a_{0}$, and thus 
$$
t_{\cT}^{-}(m , V) 
= 
\frac{\norm{m - p_{\cT}^{-}(m , V)}}{\norm{V}} 
\geq 
\frac{\a_{0}}{\norm{V}}~.
$$ 

\smallskip

In addition, since $m , p_{\cQ}^{-}(m , V) \in \clos{\cQ}$ and $\k_{0} = \diam{d}{\clos{\cQ}}$, 
we also have 
$$
t_{\cQ}^{-}(m , V) 
= 
\frac{\norm{m - p_{\cQ}^{-}(m , V)}}{\norm{V}} 
\leq 
\frac{\k_{0}}{\norm{V}}~.
$$ 

\smallskip

Hence, 
\begin{equation} \label{equ:c1-t-} 
   \frac{t_{\cQ}^{-}(m , V)}{t_{\cT}^{-}(m , V)} 
   \leq 
   \frac{\k_{0}}{\a_{0}}~. 
\end{equation} 

\medskip

Finally, if 
$K_{1} = K_{1}(a , b , c) 
\as \min{\! \{ \a_{0} / \k_{0} \ , \ \bar{m_{0} p_{0}} / \bar{m_{0} q_{0}} \}} \in (0 , 1]$, 
Equations~\ref{equ:c1-t+} and~\ref{equ:c1-t-} lead to 
$$
\frac{1}{t_{\cQ}^{-}(m , V)} + \frac{1}{t_{\cQ}^{+}(m , V)} 
\geq 
K_{1} \! \left( \frac{1}{t_{\cT}^{-}(m , V)} + \frac{1}{t_{\cT}^{+}(m , V)} \right),
$$ 
or equivalently 
$$
\FQ(m , V) \geq K_{1} \FT(m , V).
$$ 

\bigskip 

$\bullet$ \textsf{Case~2:} $V \in S(V_{2} , V_{3})$ (see Figure~\ref{fig:case-2}). 

\begin{figure}[h]
   \includegraphics[width=10cm,height=13cm,keepaspectratio=true]{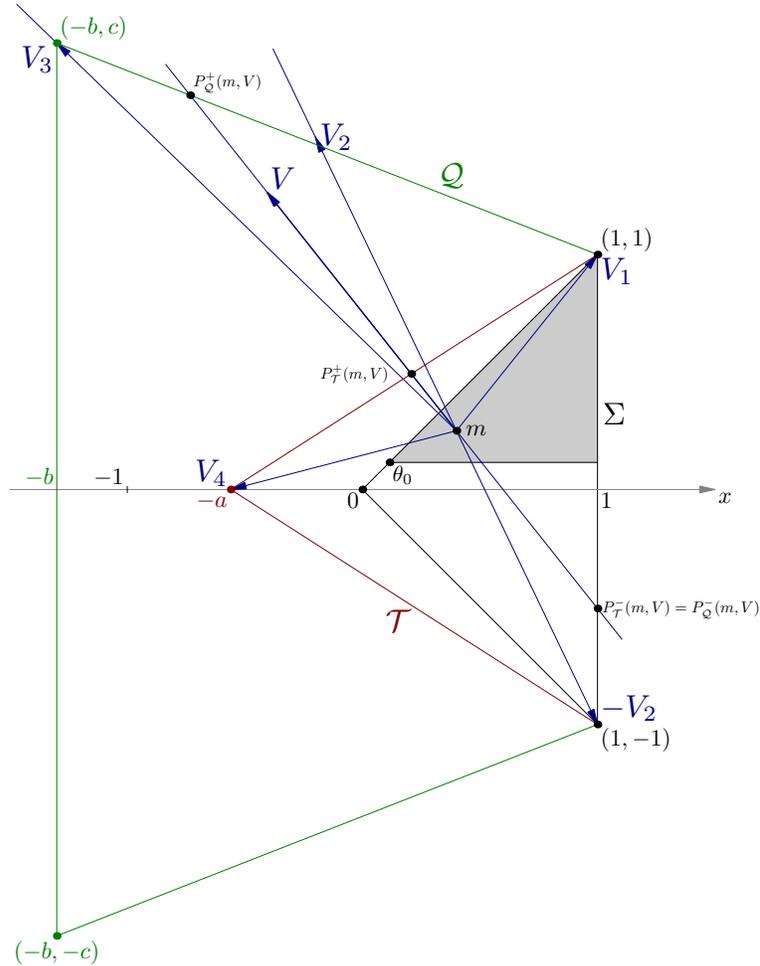}
   \caption{\label{fig:case-2} The case when $V \in S(V_{2} , V_{3})$}
\end{figure}

\bigskip

Since $V_{2} \in S(V_{1} , V_{3})$, we have $V \in S(V_{2} , V_{3}) \inc S(V_{1} , V_{3})$, 
and hence the half-line $m + \RR_{+} V$ intersects with the segment $[(1 , 1) \, , \, (-b , c)] \inc \bQ$. 

\smallskip

On the other hand, this half-line also intersects with the segment $[(1 , 1) \, , \, (-a , 0)] \inc \bT$ 
since $V \in S(V_{2} , V_{3}) \inc S(V_{1} , V_{4})$ (indeed, we have $V_{2} , V_{3} \in S(V_{1} , V_{4})$). 

\smallskip

Then, writing that the point $p_{\cT}^{+}(m , V) = (x + t_{\cT}^{+}(m , V) \l \ , \ y + t_{\cT}^{+}(m , V) \m)$ 
(resp.~$p_{\cQ}^{+}(m , V) = (x + t_{\cQ}^{+}(m , V) \l \ , \ y + t_{\cQ}^{+}(m , V) \m)$) belongs to 
the line passing through $(1 , 1)$ and $(-a , 0)$ (resp. $(-b , c)$) whose equation is 
$X - (1 + a) Y + a = 0$ (resp.~$(1 - c) X - (1 + b) Y + (b + c) = 0$), we get 
$(1 + a) \m - \l > 0$ and $(c - 1) \l + (1 + b) \m > 0$ together with 
\begin{equation} \label{equ:c2-t+} 
   t_{\cT}^{+}(m , V) 
   = 
   \frac{x - (1 + a) y + a}{(1 + a) \m - \l} 
   \quad \mbox{and} \quad 
   t_{\cQ}^{+}(m , V) 
   = 
   \frac{(1 - c) x - (1 + b) y + (b + c)}{(c - 1) \l + (1 + b) \m}~. 
\end{equation} 

\medskip

Next, since $-V_{3} \in S(-V_{2} \, , V_{1})$, we have 
$-V \in S(-V_{2} \, , -V_{3}) \inc S(-V_{2} \, , V_{1})$, 
and hence the half-line $m - \RR_{+} V$ intersects with the segment 
$[(1 , -1) \, , \, (1 , 1)] \inc \bT \cap \bQ$. 

\smallskip

The points $p_{\cT}^{-}(m , V) = (x - t_{\cT}^{-}(m , V) \l \ , \ y - t_{\cT}^{-}(m , V) \m)$ 
and $p_{\cQ}^{-}(m , V) = (x - t_{\cQ}^{-}(m , V) \l \ , \ y - t_{\cQ}^{-}(m , V) \m)$) thus lie on 
the line passing through $(1 , -1)$ and $(1 , 1)$ whose equation is 
$X - 1 = 0$, which gives $\l < 0$ and 
\begin{equation} \label{equ:c2-t-} 
   t_{\cT}^{-}(m , V) 
   = 
   t_{\cQ}^{-}(m , V) 
   = 
   \frac{x - 1}{\l}~. 
\end{equation} 

\medskip

Now, from Equations~\ref{equ:c2-t+} and~\ref{equ:c2-t-}, one obtains 
$$
2 \FQ(m , V) 
= 
\frac{1}{t_{\cQ}^{-}(m , V)} + \frac{1}{t_{\cQ}^{+}(m , V)} 
= 
\frac{1 + b}{x - 1} \times \frac{\l (1 - y) - \m (1 - x)}{(1 - c) x - (1 + b) y + (b + c)}
$$ 
together with 
$$
2 \FT(m , V) 
= 
\frac{1}{t_{\cT}^{-}(m , V)} + \frac{1}{t_{\cT}^{+}(m , V)} 
= 
\frac{1 + a}{x - 1} \times \frac{\l (1 - y) - \m (1 - x)}{x - (1 + a) y + a},
$$ 
and hence 
\begin{equation} \label{equ:c2-FQ/FC} 
   \frac{\FQ(m , V)}{\FT(m , V)} 
   = 
   \frac{1 + b}{1 + a} \times \frac{x - (1 + a) y + a}{(1 - c) x - (1 + b) y + (b + c)}~. 
\end{equation} 

\medskip

As $y \leq x < 1$, we have both $x - (1 + a) y + a \geq a (1 - y) > 0$ 
and $0 < c (1 - x) + b (1 - y) + (x - y) = (1 - c) x - (1 + b) y + (b + c) \leq (b + c) (1 - y)$ 
(indeed, $(1 - c) x \leq (1 - c) y$ since $c > 1$), from which Equation~\ref{equ:c2-FQ/FC} finally yields 
$$
\FQ(m , V) \geq K_{2} \FT(m , V),
$$ 
where $\disp K_{2} = K_{2}(a , b , c) \as \frac{a (1 + b)}{(1 + a) (b + c)} \in (0 , 1]$. 

\bigskip 

$\bullet$ \textsf{Case~3:} $V \in S(V_{3} , V_{4})$ (see Figure~\ref{fig:case-3}). 

\begin{figure}[h]
   \includegraphics[width=10cm,height=13cm,keepaspectratio=true]{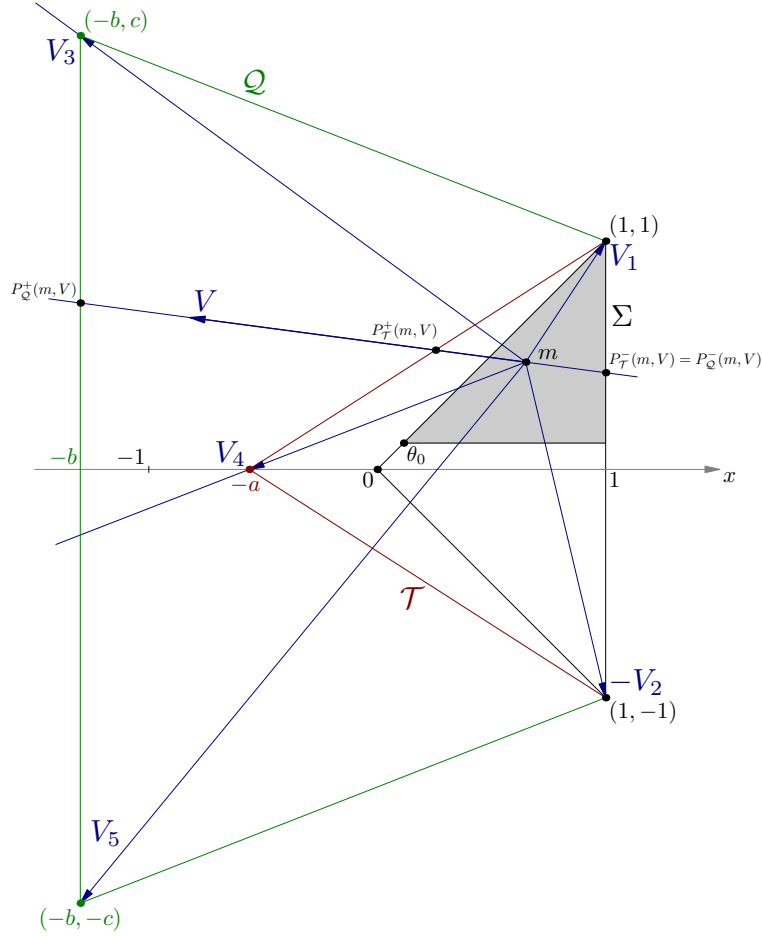}
   \caption{\label{fig:case-3} The case when $V \in S(V_{3} , V_{4})$}
\end{figure}

\bigskip

Since $-V_{3} , -V_{4} \in S(-V_{2} \, , V_{1})$, we have 
$-V \in S(-V_{3} \, , -V_{4}) \inc S(-V_{2} \, , V_{1})$, 
and hence the half-line $m - \RR_{+} V$ intersects with the segment 
$[(1 , -1) \, , \, (1 , 1)] \inc \bT \cap \bQ$. 

\smallskip

So, we get again $\l < 0$ and 
\begin{equation} \label{equ:c3-t-} 
   t_{\cT}^{-}(m , V) 
   = 
   t_{\cQ}^{-}(m , V) 
   = 
   \frac{x - 1}{\l}~. 
\end{equation} 

\medskip

On the other hand, let $V_{5} \as (-b , -c) - m = (-b - x \ , \ -c - y)$. 

\smallskip

We have $\detb{(e_{1} , e_{2})}{(V_{3} , V_{4})} \geq 0$, 
$\detb{(e_{1} , e_{2})}{(V_{4} , V_{5})} = c x + (a - b) y + a c \geq (a + c - b) y + a c \geq 0$ 
(since $x \geq y$ and $c \geq b$) and $\detb{(e_{1} , e_{2})}{(V_{3} , V_{5})} = 2 c (b + x) > 0$, 
which shows that $(V_{3} , V_{5})$ is a basis of $\Rn{2}$ having the same orientation as $(e_{1} , e_{2})$ 
with $V_{4} \in S(V_{3} , V_{5})$. 
Therefore, $V \in S(V_{3} , V_{4}) \inc S(V_{3} , V_{5})$, and hence the half-line $m + \RR_{+} V$ 
intersects with the segment $[(-b , c) \, , \, (-b , -c)] \inc \bQ$. 

\smallskip

Furthermore, this half-line also intersects with the segment $[(1 , 1) \, , \, (-a , 0)] \inc \bT$ 
since $V \in S(V_{3} , V_{4}) \inc S(V_{1} , V_{4})$ (indeed, we have $V_{3} \in S(V_{1} , V_{4})$). 

\smallskip

Then, writing that the point $p_{\cT}^{+}(m , V)$ (resp.~$p_{\cQ}^{+}(m , V)$) belongs to 
the line passing through $(1 , 1)$ and $(-a , 0)$ (resp.~$(-b , c)$ and $(-b , -c)$) 
whose equation is $X - (1 + a) Y + a = 0$ (resp.~$X + b = 0$), we compute 
$(1 + a) \m - \l > 0$ together with 
\begin{equation} \label{equ:c3-t+} 
   t_{\cT}^{+}(m , V) 
   = 
   \frac{x - (1 + a) y + a}{(1 + a) \m - \l} 
   \quad \mbox{and} \quad 
   t_{\cQ}^{+}(m , V) 
   = 
   \frac{x + b}{-\l}~. 
\end{equation} 

\medskip

Equations~\ref{equ:c3-t-} and~\ref{equ:c3-t+} then yield 
$\l (1 - y) - \m (1 - x) < 0$ and 
\begin{equation} \label{equ:c3-FQ/FC} 
   \frac{\FQ(m , V)}{\FT(m , V)} 
   = 
   \frac{1 + b}{1 + a} 
   \times 
   \frac{x - (1 + a) y + a}{x + b} 
   \times 
   \frac{\l}{\l (1 - y) - \m (1 - x)}~. 
\end{equation} 

\medskip

Since $V \in S(V_{3} , V_{4})$ and $(V_{3} , V_{4})$ is a basis of $\RR^{2}$ 
with the same orientation as $(e_{1} , e_{2})$, we have 
$\detb{(e_{1} , e_{2})}{(V_{3} , V)} = (\l y - \m x) - (c \l + b \m) \geq 0$ 
and $\detb{(e_{1} , e_{2})}{(V , V_{4})} = (\l y - \m x) - a \m \leq 0$, 
and thus $-(b - a) \m \geq c \l$. 

\smallskip

But $y \leq x < 1$, $b > a$ and $\l < 0$ then imply 
\begin{eqnarray*} 
   0 > (b - a) (\l (1 - y) - \m (1 - x)) 
   & = & 
   (b - a) \l (1 - y) - (b - a) \m (1 - x) \\ 
   & \geq & 
   (b - a) \l (1 - y) + c \l (1 - x) \\ 
   & \geq & 
   (b - a) \l (1 - y) + c \l (1 - y) = \l (b - a + c) (1 - y), 
\end{eqnarray*} 
and hence 
$$
\frac{\l}{\l (1 - y) - \m (1 - x)} 
= 
\frac{(b - a) \l}{(b - a) (\l (1 - y) - \m (1 - x))} 
\geq 
\frac{b - a}{(b - a + c) (1 - y)} > 0.
$$ 

\smallskip

Finally, using $x - (1 + a) y + a \geq a (1 - y) > 0$ together with $0 < x + b \leq 1 + b$, 
Equation~\ref{equ:c3-FQ/FC} gives 
$$
\FQ(m , V) \geq K_{3} \FT(m , V),
$$ 
where $\disp K_{3} = K_{3}(a , b , c) \as \frac{a (b - a)}{(1 + a) (b - a + c)} \in (0 , 1]$. 

\bigskip 

$\bullet$ \textsf{Case~4:} $V \in S(V_{4} , -V_{1})$ (see Figure~\ref{fig:case-4}). 

\begin{figure}[h]
   \includegraphics[width=10cm,height=13cm,keepaspectratio=true]{./Figures/polygons-fig-9.eps}
   \caption{\label{fig:case-4} The case when $V \in S(V_{4} , -V_{1})$}
\end{figure}

\bigskip

Since $-V_{4} \in S(-V_{2} \, , V_{1})$, we have $-V \in S(-V_{4} \, , V_{1}) \inc S(-V_{2} \, , V_{1})$, 
and hence the half-line $m - \RR_{+} V$ intersects with the segment 
$[(1 , -1) \, , \, (1 , 1)] \inc \bT \cap \bQ$. 

\smallskip

So, $t_{\cT}^{-}(m , V) = t_{\cQ}^{-}(m , V)$, which means that 
\begin{equation} \label{equ:c4-t-} 
   \frac{t_{\cQ}^{-}(m , V)}{t_{\cT}^{-}(m , V)} = 1. 
\end{equation} 

\medskip

On the other hand, since $m \in \S$ 
and $p_{\cT}^{+}(m , V) \in [(-a , 0) \, , \, (1 , -1)] \inc \RR \cart (-\infty , 0]$ 
(indeed, $-V_{1} \in S(V_{4} , -V_{2})$ implies $V \in S(V_{4} , -V_{1}) \inc S(V_{4} , -V_{2})$), 
we have $\norm{m - p_{\cT}^{+}(m , V)} \geq d(m , \RR \cart \{ 0 \}) \geq \a_{0}$, and thus 
$$
t_{\cT}^{+}(m , V) 
= 
\frac{\norm{m - p_{\cT}^{+}(m , V)}}{\norm{V}} 
\geq 
\frac{\a_{0}}{\norm{V}}~.
$$ 
In addition, since $m , p_{\cQ}^{+}(m , V) \in \clos{\cQ}$ and $\k_{0} = \diam{d}{\clos{\cQ}}$, 
we also have 
$$
t_{\cQ}^{+}(m , V) 
= 
\frac{\norm{m - p_{\cQ}^{+}(m , V)}}{\norm{V}} 
\leq 
\frac{\k_{0}}{\norm{V}}~.
$$ 

\smallskip

Hence, 
\begin{equation} \label{equ:c4-t+} 
   \frac{t_{\cQ}^{+}(m , V)}{t_{\cT}^{+}(m , V)} 
   \leq 
   \frac{\k_{0}}{\a_{0}}~. 
\end{equation} 

\medskip

Finally, if 
$K_{4} = K_{4}(a , b , c) \as \min{\! \{ 1 \ , \ \a_{0} / \k_{0} \}} \in (0 , 1]$, 
Equations~\ref{equ:c4-t-} and~\ref{equ:c4-t+} lead to 
$$
\FQ(m , V) \geq K_{4} \FT(m , V).
$$ 

\bigskip 

At this stage of the proof, defining 
$K = K(a , b , c) \as \min{\! \{ K_{1} , K_{2} , K_{3} , K_{4} \}} \in (0 , 1]$ 
and summing up the results obtained in the four cases discussed above, we can write 
\begin{equation} \label{equ:part-1} 
   \FQ(m , V) \geq K \FT(m , V) 
\end{equation} 
for all $m \in \S$ and $V \in T_{\! m}\cT = T_{\! m}\cQ = \Rn{2}$. 

\bigskip

Now, the only thing to be done is to establish a similar inequality as in Equation~\ref{equ:part-1} 
for $m \in \D^{\! +} \setmin \S$, from which we will get Proposition~\ref{prop:comparison-1} 
since both $\cT$ and $\cQ$ are preserved by the reflection about the $x$-axis. 

\smallskip

So, let $\d_{0} \as d(\D^{\! +} \setmin \S \ , \ [(1 , 1) \, , \, (-a , 0)]) > 0$ 
and consider a point $m \in \D^{\! +} \setmin \S$ together with a vector $V = (\l , \m) \in \Rn{2}$ 
such that $V \neq 0$. 

\smallskip

First of all, since $m , p_{\cQ}^{+}(m , V) \in \clos{\cQ}$ and $\k_{0} = \diam{d}{\clos{\cQ}}$, 
we have 
$$
t_{\cQ}^{+}(m , V) 
= 
\frac{\norm{m - p_{\cQ}^{+}(m , V)}}{\norm{V}} 
\leq 
\frac{\k_{0}}{\norm{V}}~.
$$ 

\smallskip

Next, the Finsler metrics $\FT$ and $\FQ$ being reversible, we can assume that $\l \leq 0$, 
and hence $p_{\cT}^{+}(m , V) \in \bT \setmin (\{ 1 \} \cart [-1 , 1])$. 

\smallskip 

This implies that $\norm{m - p_{\cT}^{+}(m , V)} \geq \d_{0}$ and gives 
$$
t_{\cT}^{+}(m , V) 
= 
\frac{\norm{m - p_{\cT}^{+}(m , V)}}{\norm{V}} 
\geq 
\frac{\d_{0}}{\norm{V}}~.
$$ 

\smallskip

Therefore, 
\begin{equation} \label{equ:t+} 
   \frac{t_{\cQ}^{+}(m , V)}{t_{\cT}^{+}(m , V)} 
   \leq 
   \frac{\k_{0}}{\d_{0}}~. 
\end{equation} 

\medskip

On the other hand, as regards $t_{\cT}^{-}(m , V)$ and $t_{\cQ}^{-}(m , V)$, we have two cases to look at. 

\bigskip

$\bullet$ \textsf{First case:} $p_{\cT}^{-}(m , V) \in \{ 1 \} \cart [-1 , 1] \inc \bT$. 

\smallskip

In that case, we also have $p_{\cQ}^{-}(m , V) \in \{ 1 \} \cart [-1 , 1] \inc \bQ$, 
and hence $t_{\cT}^{-}(m , V) = t_{\cQ}^{-}(m , V)$, or equivalently 
\begin{equation} \label{equ:1-t-} 
   \frac{t_{\cQ}^{-}(m , V)}{t_{\cT}^{-}(m , V)} = 1. 
\end{equation} 

\medskip

So, if $K_{5} = K_{5}(a , b , c) \as \min{\! \{ 1 \ , \ \d_{0} / \k_{0} \}} \in (0 , 1]$, 
Equations~\ref{equ:t+} and~\ref{equ:1-t-} yield 
\begin{equation} \label{equ:part-2-1} 
   \FQ(m , V) \geq K_{5} \FT(m , V). 
\end{equation} 

\bigskip

$\bullet$ \textsf{Second case:} $p_{\cT}^{-}(m , V) \in \bT \setmin (\{ 1 \} \cart [-1 , 1])$. 

\smallskip

Then we have 
$$
t_{\cT}^{-}(m , V) 
= 
\frac{\norm{m - p_{\cT}^{-}(m , V)}}{\norm{V}} 
\geq 
\frac{\d_{0}}{\norm{V}}~.
$$ 

\smallskip

But 
$$
t_{\cQ}^{-}(m , V) 
= 
\frac{\norm{m - p_{\cQ}^{-}(m , V)}}{\norm{V}} 
\leq 
\frac{\k_{0}}{\norm{V}}
$$ 
since $m , p_{\cQ}^{+}(m , V) \in \clos{\cQ}$ and $\k_{0} = \diam{d}{\clos{\cQ}}$. 

\smallskip

Therefore, 
\begin{equation} \label{equ:2-t-} 
   \frac{t_{\cQ}^{-}(m , V)}{t_{\cT}^{-}(m , V)} 
   \leq 
   \frac{\k_{0}}{\d_{0}}~. 
\end{equation} 

\medskip

If $K_{6} = K_{6}(a , b , c) \as \d_{0} / \k_{0} \in (0 , 1]$, 
Equations~\ref{equ:t+} and~\ref{equ:2-t-} thus lead to 
\begin{equation} \label{equ:part-2-2} 
   \FQ(m , V) \geq K_{6} \FT(m , V). 
\end{equation} 

\bigskip

Conclusion: combining Equations~\ref{equ:part-1}, \ref{equ:part-2-1} and~\ref{equ:part-2-2}, 
and defining 
$$
A = A(a , b , c) \as \min{\! \{ K , K_{5} , K_{6} \}} \in (0 , 1],
$$ 
we have finally obtained that 
$$
\FQ(m , V) \geq A \FT(m , V)
$$ 
for all $m \in \D^{\! +}$ and $V \in T_{\! m}\cT = T_{\! m}\cQ = \Rn{2}$, 
which ends the proof of Proposition~\ref{prop:comparison-1}. 
\end{proof} 

\bigskip

From Proposition~\ref{prop:comparison-1}, we can then deduce 

\begin{proposition} \label{prop:comparison-2} 
   Let $\cC_{1 \!\!}$ and $\cC_{2}$ be open bounded convex sets in $\Rn{2}$ such that 
   
   \begin{enumerate}
      \item the segment $\{ 1 \} \cart [-1 , 1]$ 
      is included in both boundaries $\bC_{1 \!\!}$ and $\bC_{2}$, 
      
      \smallskip
      
      \item $(1 , -1)$ and $(1 , 1)$ are corner points of\;\;{}{}$\clos{\cC}_{\! 1 \!}$ 
      and\;\;{}$\clos{\cC}_{\! 2}$, 
      
      \smallskip
      
      \item the origin $0$ lies in $\cC_{1} \cap \cC_{2}$, and 
      
      \smallskip
      
      \item $\D \inc \cC_{1} \cap \cC_{2}$. 
   \end{enumerate} 
   
   Then there exists a constant $B = B(\cC_{1} , \cC_{2}) \geq 1$ that satisfies 
   $$
   \frac{1}{B} F_{\cC_{1}}(m , V) \leq F_{\cC_{2}}(m , V) \leq B F_{\cC_{1}}(m , V)
   $$ 
   for all $m \in \D$ and $V \in T_{\! m}\cC_{1} = T_{\! m}\cC_{2} = \Rn{2}$. 
\end{proposition}

\medskip

\begin{proof}~\\ 
Since $\cC_{1} \cap \cC_{2}$ is an open set in $\Rn{2}$ 
that contains the origin $0$ by point~(3), its intersection $I$ with $\RR \cart \{ 0 \}$ 
is an open set in $\RR \cart \{ 0 \}$ which also contains $0$, 
and hence there exists a number $a \in (0 , 1)$ such that $[-a , a] \cart \{ 0 \} \inc I$. 
This implies that $(-a , 0) \in I \inc \cC_{1} \cap \cC_{2}$, 
and therefore $\cT \inc \cC_{1} \cap \cC_{2}$, where $\cT \inc \Rn{2}$ is the triangle defined 
as the open convex hull of the points $(1 , -1)$, $(1 , 1)$ and $(-a , 0)$ 
(indeed, $\cC_{1}$ and $\cC_{2}$ are convex sets in $\Rn{2}$ whose boundaries contain 
$(1 , -1)$ and $(1 , 1)$ by point~(1)). 

\medskip

On the other hand, since the sets $\cC_{1}$ and $\cC_{2}$ are bounded, there exists a number $b \geq 1$ 
such that they are \emph{both} inside $[-b , b] \cart [-b , b]$. 
So, $\cC_{1}$ and $\cC_{2}$ are included in the open half-plane of $\Rn{2}$ 
whose boundary is the line $\{ -b \} \cart \RR$ and which contains the origin $0$. 

\medskip

Next, the convexity of $\cC_{1}$ (resp.~$\cC_{2}$) together with points~(1) and~(3) show that 
$\cC_{1}$ (resp.~$\cC_{2}$) lies inside the open half-plane of $\Rn{2}$ 
whose boundary is the line $\RR ((1 , -1) - (1 , 1)) = \{ 1 \} \cart \RR$ 
and which contains the origin $0$. 

Moreover, point~(2) implies that $\cC_{1}$ (resp.~$\cC_{2}$) has support lines 
$\cL_{1}^{-}$ (resp.~$\cL_{2}^{-}$) and $\cL_{1}^{+}$ (resp.~$\cL_{2}^{+}$) 
at respectively $(1 , -1)$ and $(1 , 1)$ which are not equal to the line $\{ 1 \} \cart \RR$. 
Therefore, $\cC_{1}$ (resp.~$\cC_{2}$) lies inside the open half-planes of $\Rn{2}$ 
which contain the origin $0$ and whose boundaries are the lines 
$\cL_{1}^{-}$ (resp.~$\cL_{2}^{-}$) and $\cL_{1}^{+}$ (resp.~$\cL_{2}^{+}$). 

So, if we denote by $c_{1}$ (resp.~$c_{2}$) the maximum of the absolute values of the second coordinates 
of the intersection points of the lines $\cL_{1}^{-}$ (resp.~$\cL_{2}^{-}$) 
and $\cL_{1}^{+}$ (resp.~$\cL_{2}^{+}$) with the line $\{ -b \} \cart \RR$, 
then $\cC_{1}$ and $\cC_{2}$ are included in the open half-planes of $\Rn{2}$ 
which contain the origin $0$ and whose boundaries are the lines $\RR ((1 , -1) - (-b , -c))$ 
and $\RR ((-b , c) - (1 , 1))$, where $c \as \max{\! \{ c_{1} , c_{2} \}} + b + 1 > b$. 

\medskip

Conclusion: we have $\cC_{1} \inc \cQ$ and $\cC_{2} \inc \cQ$, 
where $\cQ \inc \Rn{2}$ is the quadrilateral defined as the open convex hull 
of the points $(1 , -1)$, $(1 , 1)$, $(-b , c)$ and $(-b , -c)$. 

\bigskip

Now, for all $m \in \D$ and $V \in T_{\! m}\cT = T_{\! m}\cC_{1} = T_{\! m}\cC_{2} = T_{\! m}\cQ = \Rn{2}$, 
we can write 
\begin{eqnarray*} 
   A F_{\cC_{1}}(m , V) 
   & \leq & 
   A \FT(m , V) \quad \mbox{(since $\cT \inc \cC_{1}$)} \\ 
   & \leq & 
   \FQ(m , V) \quad \mbox{(by the first inequality in Proposition~\ref{prop:comparison-1})} \\ 
   & \leq & 
   F_{\cC_{2}}(m , V) \quad \mbox{(since $\cC_{2} \inc \cQ$)} \\ 
   & \leq & 
   \FT(m , V) \quad \mbox{(since $\cT \inc \cC_{2}$)} \\ 
   & \leq & 
   \frac{1}{A} \FQ(m , V) \quad 
   \mbox{(by the first inequality in Proposition~\ref{prop:comparison-1})} \\ 
   & \leq & 
   \frac{1}{A} F_{\cC_{1}}(m , V) \quad \mbox{(since $\cC_{1} \inc \cQ$)}, 
\end{eqnarray*} 
which proves Proposition~\ref{prop:comparison-2} with $B \as 1 / A \geq 1$. 
\end{proof}

\bigskip
\bigskip
\bigskip


\section{Lipschitz equivalence to Euclidean plane} \label{sec:Lipschitz-eq} 

In this section, we build a homeomorphism from an open convex polygonal set to Euclidean plane, 
and prove that it is bi-Lipschitz with respect to the Hilbert metric of the polygonal set 
and the Euclidean distance of the plane. This is the statement of Theorem~\ref{thm:bi-Lipschitz}. 

\bigskip

So, let $\cP$ be an open convex polygonal set in $\Rn{2}$ that contains the origin $0$. 

\bigskip

Let $\llist{v}{1}{n}$ be the vertices of $\cP$ ($i.e.$, the corner points of the convex set $\cP$) 
that we assume to be cyclically ordered in $\bP$ (notice we have $n \geq 3$). 

\bigskip

Define $v_{0} \as v_{n}$ and $v_{n + 1} \as v_{1}$. 

\bigskip

Let $f : \cP \to \Rn{2}$ be the map defined as follows. 

\smallskip

For each $k \in \{ 1 , \ldots , n \}$, let 
$\D_{k} \as \{ s v_{k} + t v_{k + 1} \st s \geq 0, \ t \geq 0 \ \ \mbox{and} \ \ s + t < 1 \} \inc \cP$, 
and consider the unique linear transformation $L_{k}$ of $\Rn{2}$ 
such that $L_{k}(v_{k}) \as (1 , -1)$ and $L_{k}(v_{k + 1}) \as (1 , 1)$. 

\smallskip

Then, given any $p \in \D_{k}$, we define 
$$
f(p) \as L_{k}^{-1}(\F(L_{k}(p))),
$$ 
where $\F : \cS \to \Rn{2}$ is the map considered in Proposition~\ref{prop:map}. 

\smallskip

\begin{figure}[h]
   \includegraphics[width=15cm,height=15cm,keepaspectratio=true]{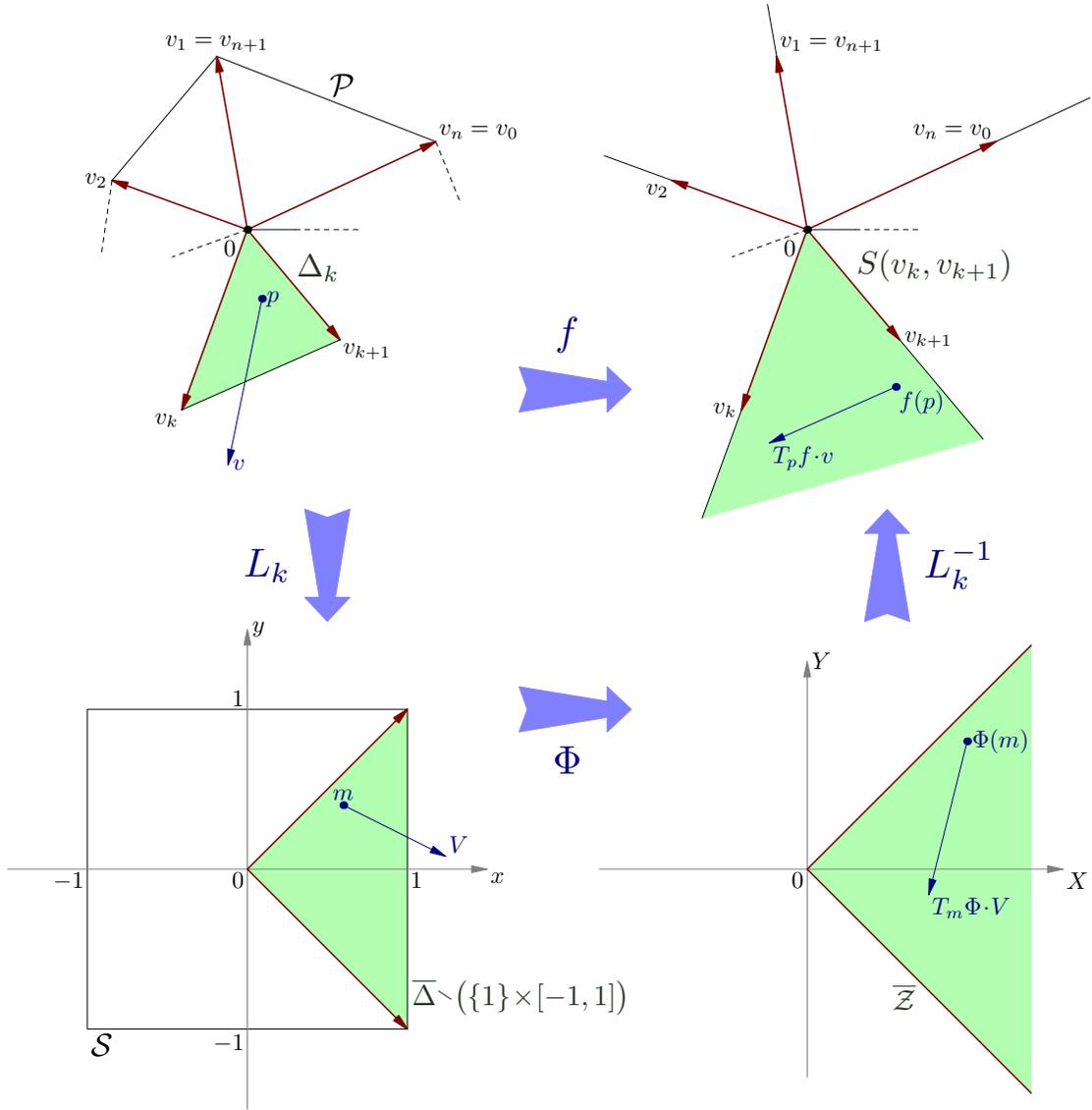}
   \caption{\label{fig:bi-Lipschitz} The bi-Lipschitz homeomorphism $f$}
\end{figure}

In other words, $f$ makes the following diagram commute (see Figure~\ref{fig:bi-Lipschitz}): 
$$
\begin{CD} 
   \D_{k} @>{f}>> S(v_{k} , v_{k + 1}) \\ 
   @V{L_{k}}VV @VV{L_{k}}V \\ 
   \clos{\D} \setmin (\{ 1 \} \cart [-1 , 1]) @>>{\F}> \clos{Z} 
\end{CD},
$$ 
where we recall that $\cS \as ]{-1} , 1[ \: \cart \: ]{-1} , 1[ \; \inc \Rn{2}$, 
$\D \as \{ (x , y) \in \Rn{2} \st |y| < x < 1 \} \inc \cS$ 
and $\cZ \as \{ (X , Y) \in \Rn{2} \st |Y| < X \} \inc \Rn{2}$. 

\medskip

This makes sense since $\disp \bigcup_{i = 1}^{n} \D_{i} = \cP$, 
$L_{k}(\D_{k}) = \clos{\D} \setmin (\{ 1 \} \cart [-1 , 1]) \inc \cS$ 
and for all $\disp p \in [0 , 1) v_{k} = \D_{k - 1} \, \cap \, \D_{k}$, 
we have $L_{k - 1}^{-1}(\F(L_{k - 1}(p))) = L_{k}^{-1}(\F(L_{k}(p))) \in \RR v_{k}$. 

\bigskip

With this definition, keeping in mind that $\norm{\cdot}$ and $d$ denote respectively 
the usual $\ll{1}$-norm on $\Rn{2}$ and its associated distance, we get 

\begin{theorem} \label{thm:bi-Lipschitz} 
   The map $f$ satisfies the following properties: 
   
   \begin{enumerate}
      \item $f$ is a homeomorphism. 
      
      \medskip
      
      \item If \! $\cU$~and~$\cV$~are the open sets in $\Rn{2}$ defined by 
      $$
      \cU \as \cP \, \setmin \! \left( \bigcup_{k = 1}^{n} [0 , 1) v_{k} \!\! \right) 
      \quad \mbox{and} \quad 
      \cV \as \Rn{2} \, \setmin \! \left( \bigcup_{k = 1}^{n} \RR v_{k} \!\! \right) \!\! ,
      $$ 
      then $f(\cU) = \cV$ and $f$ induces a smooth diffeomorphism from $\cU$ onto $\cV$. 
      
      \medskip
      
      \item There exists a constant $C \geq 1$ such that 
      $$
      \frac{1}{C} \dP(p , q) \leq d(f(p) , f(q)) \leq C \dP(p , q)
      $$ 
      for all $p , q \in \cP$. 
   \end{enumerate}
\end{theorem}

\bigskip

Before proving this theorem, let us establish the following: 

\begin{lemma} \label{lem:comparison-3} 
   Given any real number $\a > 0$, there is a constant $M = M(\a) \geq 1$ such that 
   $$
   \frac{1}{M} \ln{\!\! \left( \! \frac{1 - s}{1 - t} \times \frac{1 + \a t}{1 + \a s} \! \right)} 
   \leq 
   \ln{\!\! \left( \! \frac{1 - s}{1 - t} \times \frac{1 + t}{1 + s} \! \right)} 
   \leq 
   M \ln{\!\! \left( \! \frac{1 - s}{1 - t} \times \frac{1 + \a t}{1 + \a s} \! \right)}
   $$ 
   for all $0 \leq s < t < 1$. 
\end{lemma}

\medskip

\begin{proof}~\\ 
Consider $D \as \{ (s , t) \in \Rn{2} \st 0 \leq s < t \leq 1 \}$ 
and let $\f : D \to \RR$ be the function defined by 
$$
\f(s , t) 
\as 
\ln{\!\! \left( \! \frac{1 + t}{1 + s} \! \right)} 
\! {\bigg /} \! 
\ln{\!\! \left( \! \frac{1 + \a t}{1 + \a s} \! \right)}.
$$ 

\medskip

Given $\l \in [0 , 1]$, we have for all $(s , t) \in D$, 
$$
\f(s , t) 
\thicksim 
\! \left( \! \frac{1 + t}{1 + s} - 1 \! \right) 
\!\! {\bigg /} \!\! 
\left( \! \frac{1 + \a t}{1 + \a s} - 1 \! \right) 
= 
\frac{1 + \a s}{\a (1 + s)} \to \frac{1 + \a \l}{\a (1 + \l)}
$$ 
as $(s , t) \to (\l , \l)$. 

\medskip

Hence, by continuity of $\f$, the function $\hat{\f} : \clos{D} \to \RR$ defined by 
$$
\hat{\f}(s , t) \as \f(s , t) \ \ \mbox{if} \ \ (s , t) \in D 
\quad \mbox{and} \quad 
\disp \hat{\f}(s , t) \as \frac{1 + \a s}{\a (1 + s)} \ \ \mbox{if} \ \ s = t
$$ 
is continuous. 

Then, compactness of $\clos{D}$ implies that $\hat{\f}$ has a minimum and a maximum. 
But these latters are positive since one can easily check that $\hat{\f}(s , t) > 0$ 
for all $(s , t) \in \clos{D}$, and this implies that there is a constant $M \geq 1$ 
such that 
\begin{equation} \label{equ:bounded-map} 
   \frac{1}{M} \ln{\!\! \left( \! \frac{1 + \a t}{1 + \a s} \! \right)} 
   \leq 
   \ln{\!\! \left( \! \frac{1 + t}{1 + s} \! \right)} 
   \leq 
   M \ln{\!\! \left( \! \frac{1 + \a t}{1 + \a s} \! \right)} 
\end{equation} 
for all $(s , t) \in D$. 

\medskip

Finally, for all $0 \leq s < t < 1$, we have 
\begin{eqnarray*} 
   \frac{1}{M} \ln{\!\! \left( \! \frac{1 - s}{1 - t} \times \frac{1 + \a t}{1 + \a s} \! \right)} 
   & = & 
   \frac{1}{M} \ln{\!\! \left( \! \frac{1 - s}{1 - t} \! \right)} 
   + \frac{1}{M} \ln{\!\! \left( \! \frac{1 + \a t}{1 + \a s} \! \right)} \\ 
   & \leq & 
   \ln{\!\! \left( \! \frac{1 - s}{1 - t} \! \right)} 
   + \frac{1}{M} \ln{\!\! \left( \! \frac{1 + \a t}{1 + \a s} \! \right)} 
   \quad \mbox{(since $1 / M \leq 1$)} \\ 
   & \leq & 
   \ln{\!\! \left( \! \frac{1 - s}{1 - t} \! \right)} 
   + \ln{\!\! \left( \! \frac{1 + t}{1 + s} \! \right)} 
   = \ln{\!\! \left( \! \frac{1 - s}{1 - t} \times \frac{1 + t}{1 + s} \! \right)} \\ 
   & & 
   \mbox{(by the first inequality in Equation~\ref{equ:bounded-map})} \\ 
   & \leq & 
   \ln{\!\! \left( \! \frac{1 - s}{1 - t} \! \right)} 
   + M \ln{\!\! \left( \! \frac{1 + \a t}{1 + \a s} \! \right)} \\ 
   & & 
   \mbox{(by the second inequality in Equation~\ref{equ:bounded-map})} \\ 
   & \leq & 
   M \ln{\!\! \left( \! \frac{1 - s}{1 - t} \! \right)} 
   + M \ln{\!\! \left( \! \frac{1 + \a t}{1 + \a s} \! \right)} 
   \quad \mbox{(since $M \geq 1$)} \\ 
   & = & 
   M \ln{\!\! \left( \! \frac{1 - s}{1 - t} \times \frac{1 + \a t}{1 + \a s} \! \right)} 
\end{eqnarray*} 

\medskip

This proves Lemma~\ref{lem:comparison-3}. 
\end{proof}

\bigskip

\begin{proof}[Proof of Theorem~\ref{thm:bi-Lipschitz}]~\\ 
$\bullet$ \textsf{Point~(1):} 
Let $g : \Rn{2} \to \cP$ be the map given by $g(P) \as L_{k}^{-1}(\F^{-1}(L_{k}(P)))$ 
for all $k \in \{ 1 , \ldots , n \}$ and $P \in S(v_{k} , v_{k + 1})$, 
this definition making sense since $\disp \bigcup_{i = 1}^{n} S(v_{i} , v_{i + 1}) = \Rn{2}$, \linebreak 
$\F^{-1}(L_{k}(S(v_{k} , v_{k + 1}))) = \clos{\D} \setmin (\{ 1 \} \cart [-1 , 1])$ 
(using $L_{k}(S(v_{k} , v_{k + 1})) = \clos{\cZ}$ and point~(1) in Proposition~\ref{prop:map}) 
and $L_{k - 1}^{-1}(\F^{-1}(L_{k - 1}(P))) = L_{k}^{-1}(\F^{-1}(L_{k}(P))) \in [0 , 1) v_{k}$ 
whenever $P \in \RR v_{k} = S(v_{k - 1} , v_{k}) \, \cap \, S(v_{k} , v_{k + 1})$. 

Then it is easy to check that $f \circ g = \I{\mathbf{R}^{\! 2}}$ and $g \circ f = \I{\cP}$ (identity maps), 
which shows that $f$ is bijective with $f^{-1} = g$. 

\medskip

In addition, $f$ and $g$ are continuous since $\llist{L}{1}{n}$ and $\F$ are homeomorphisms. 

\bigskip

$\bullet$ \textsf{Point~(2):} 
For each $k \in \{ 1 , \ldots , n \}$, we have $f([0 , 1) v_{k}) = \RR v_{k}$, 
and therefore 
$$
f(\cU) 
= 
f(\cP) \, \setmin \! \left( \bigcup_{k = 1}^{n} f([0 , 1) v_{k}) \!\! \right) 
= 
\Rn{2} \, \setmin \! \left( \bigcup_{k = 1}^{n} \RR v_{k} \!\! \right) 
= 
\cV.
$$ 

Moreover, since 
$$
\cU = \bigcup_{k = 1}^{n} \intr{\D}_{k} 
\quad \mbox{and} \quad 
\cV = \bigcup_{k = 1}^{n} \overbrace{S(v_{k} , v_{k + 1})}^{\circ}
$$ 
together with the fact that $\llist{L}{1}{n}$ are smooth by linearity 
and $\F$ is a smooth diffeomorphism by Proposition~\ref{prop:map}, 
we get that $\rest{f}{\cU}$ and $\rest{g}{\cV}$ are smooth. 

\bigskip

$\bullet$ \textsf{Point~(3):} 
Fixing $k \in \{ 1 , \ldots , n \}$ and applying Proposition~\ref{prop:comparison-2} 
with $\cC_{1} \as \cS$ and $\cC_{2} \as L_{k}(\cP)$, 
we get a constant $B_{k} \geq 1$ such that for all $m \in \D$ 
and $V \in T_{\! m}\cS = T_{\! m}(L_{k}(\cP)) = \Rn{2}$, 
$$
\frac{1}{B_{k}} \FS(m , V) 
\leq 
F_{L_{k}(\cP)}(m , V) 
\leq 
B_{k} \FS(m , V),
$$ 
and hence 
$$
\frac{1}{B_{k}} F_{L_{k}(\cP)}(m , V) 
\leq 
\norm{\Tg{m}{\F}{V}} 
\leq 
2 B_{k} F_{L_{k}(\cP)}(m , V)
$$ 
by Proposition~\ref{prop:map}. 

But, since $L_{k}$ induces an isometry from $(\cP , \dP)$ onto $(\Rn{2} , d)$ 
(being affine, $L_{k}$ preserves the cross ratio), 
this is equivalent to saying that for all $p \in \intr{\D}_{k}$ 
and $v \in T_{\! p}\cP = \Rn{2}$ (writing $m = L_{k}(p)$ and $V = L_{k}(v)$), we have 
$$
\frac{1}{B_{k}} \FP(p , v) 
\leq 
\norm{\Tg{L_{k}(p)}{\F}{L_{k}(v)}} 
\leq 
2 B_{k} \FP(p , v),
$$ 
which yields 
$$
\frac{1}{B_{k} \opnorm{L_{k}}} \FP(p , v) 
\leq 
\norm{L_{k}^{-1}(\Tg{L_{k}(p)}{\F}{L_{k}(v)})} = \norm{\Tg{p}{f}{v}} 
\leq 
2 B_{k} \opnorm{L_{k}^{-1}} \FP(p , v),
$$ 
where $\opnorm{\cdot}$ denotes the operator norm on $\End{\Rn{2}}$ associated with $\norm{\cdot}$. 

\medskip

Now, if $K \as \max{\! \left\{ B_{k} \opnorm{L_{k}} + 2 B_{k} \opnorm{L_{k}^{-1}} + 1 
~{\big |}~ 1 \leq k \leq n \right\}} \geq 1$, 
then 
\begin{equation} \label{equ:Finsler} 
   \frac{1}{K} \FP(p , v) 
   \leq 
   \norm{\Tg{p}{f}{v}} 
   \leq 
   K \FP(p , v) 
\end{equation} 
for all $\disp p \in \bigcup_{k = 1}^{n} \intr{\D}_{k} = \cU$ and $v \in T_{\! p}\cP = \Rn{2}$. 

\medskip

We will then prove Theorem~\ref{thm:bi-Lipschitz} using the fact that $(\cP , \dP)$ and $(\Rn{2} , d)$ 
are geodesic metric spaces in which affine segments are geodesics (see Introduction). 

\medskip

Let $p , q \in \cP$ and $\c : [0 , 1] \to \cP$ defined by $\c(t) \as  (1 - t) p + t q$. 
Assume that $\rest{\c}{[0 , 1)} \inc \cU$ and $\disp q = \c(1) \in \bigcup_{k = 1}^{n} [0 , 1) v_{k}$. 
The second inequality in Equation~\ref{equ:Finsler} then implies that for all $t \in [0 , 1)$, we have 
$$
d(f(p) , f(\c(t))) 
\leq 
\int{0}{t}{\norm{\Tg{\c(s)}{f}{\c'(s)}} \!}{s} 
\leq 
K \!\! \int{0}{t}{\FP(\c(s) , \c'(s))}{s} 
= 
K \dP(p , \c(t)),
$$ 
and thus $d(f(p) , f(\c(t))) \leq K \dP(p , \c(t))$, 
which gives 
\begin{equation} \label{equ:p-q} 
   d(f(p) , f(q)) \leq K \dP(p , q) 
\end{equation} 
as $t \to 1$. 

\medskip

If now $p \as s v_{k}$ and $q \as t v_{k}$ for some $0 \leq s < t < 1$ and $k \in \{ 1 , \ldots , n \}$, 
a straightforward calculation gives 
\begin{eqnarray*} 
   f(p) = L_{k}^{-1}(\F(s L_{k}(v_{k}))) 
   & = & 
   L_{k}^{-1}(\F(s (1 , -1))) \\ 
   & = & 
   L_{k}^{-1}(\F(s , -s)) = L_{k}^{-1}(\atanh{\! (s)} , -\atanh{\! (s)}) \\ 
   & = & 
   \atanh{\! (s)} L_{k}^{-1}(1 , -1) = \atanh{\! (s)} v_{k}, 
\end{eqnarray*} 
and hence 
\begin{equation} \label{equ:distance-d} 
   d(f(p) , f(q)) = \norm{f(q) - f(p)} 
   = (\atanh{\! (t)} - \atanh{\! (s)}) \norm{v_{k}} 
   = \frac{\norm{v_{k}}}{2} \ln{\!\! \left( \! \frac{1 - s}{1 - t} \times \frac{1 + t}{1 + s} \! \right)}, 
\end{equation} 
together with 
$$
e^{2 \dP(p , q)} = \left[ v_{k} \ , \ q \ , \ p \ , \ -t_{\cP}^{-}(0 , v_{k}) v_{k} \right] 
= \frac{1 - s}{1 - t} \times \frac{t + t_{\cP}^{-}(0 , v_{k})}{s + t_{\cP}^{-}(0 , v_{k})}~,
$$ 
or equivalently 
\begin{equation} \label{equ:distance-dP} 
   \dP(p , q) 
   = 
   \frac{1}{2} \ln{\!\! \left( \! \frac{1 - s}{1 - t} \times \frac{1 + t \a_{k}}{1 + s \a_{k}} \! \right)}, 
\end{equation} 
where $\a_{k} \as 1 / t_{\cP}^{-}(0 , v_{k}) > 0$. 

\medskip

Then, using Lemma~\ref{lem:comparison-3} with $\a \as \a_{k}$ and denoting 
$\L_{k} \as M \! (\a_{k}) \times \max{\! \left\{ \norm{v_{k}} \ , \ 1 / \! \norm{v_{k}} \right\}} \geq 1$, 
Equations~\ref{equ:distance-d} and~\ref{equ:distance-dP} yield 
\begin{equation} \label{equ:0-p} 
   \frac{1}{\L_{k}} \dP(p , q) \leq d(f(p) , f(q)) \leq \L_{k} \dP(p , q). 
\end{equation}

\medskip

If now $p$ and $q$ are arbitrary chosen in $\cP$, the closed affine segment joining $p$ and $q$ 
either meets $\disp \bigcup_{k = 1}^{n} [0 , 1) v_{k}$ in \emph{at most} $n$ points, 
or has at least two distinct points in common with some $[0 , 1) v_{k}$ for $k \in \{ 1 , \ldots , n \}$. 

Therefore, it follows from Equation~\ref{equ:p-q} and the second inequality in Equation~\ref{equ:0-p} 
that \linebreak 
$d(f(p) , f(q)) \leq C \dP(p , q)$ holds 
with $C \as \max{\! \{ \L_{k} \st 1 \leq k \leq n \}} + K \geq 1$. 

\medskip

Finally, using the first inequalities in Equations~\ref{equ:Finsler} and~\ref{equ:0-p}, 
the same arguments for $f^{-1}$ as those for $f$ lead to 
$\disp \dP(p , q) \leq C d(f(p) , f(q))$ for any $p , q \in \cP$. 

\medskip

This ends the proof of Theorem~\ref{thm:bi-Lipschitz}. 
\end{proof}

\bigskip
\bigskip
\bigskip


\bibliographystyle{acm}
\bibliography{./math-biblio}

\end{document}

%% file: latex-macros.tex
\theoremstyle{plain}
\newtheorem{theorem}{Theorem}[section]
\newtheorem*{theorem*}{Theorem} 

\newtheorem{lemma}{Lemma}[section]
\newtheorem*{lemma*}{Lemma} 

\newtheorem*{corollary*}{Corollary} 

\newtheorem*{consequence*}{Consequence} 
\newtheorem{proposition}{Proposition}[section]
\newtheorem*{proposition*}{Proposition} 

\newtheorem*{conjecture*}{Conjecture} 

\theoremstyle{definition}
\newtheorem{definition}{Definition}[section]
\newtheorem*{definition*}{Definition} 

\newtheorem*{remark*}{Remark} 
\newtheorem*{remarks*}{Remarks} 
\newtheorem*{notations}{Notations} 

\newtheorem*{question*}{Question} 
\newtheorem*{questions*}{Questions} 

\newtheorem*{example*}{Example} 
\newtheorem*{examples*}{Examples} 

\newtheorem*{exercise*}{Exercise} 
\newtheorem*{exercises*}{Exercises} 

\theoremstyle{plain}

\newtheorem*{thm*}{Théorème} 

\newtheorem*{lemme*}{Lemme} 

\newtheorem*{cor*}{Corollaire} 

\newtheorem*{csq*}{Conséquence} 

\newtheorem*{prop*}{Proposition} 

\newtheorem*{conj*}{Conjecture} 

\theoremstyle{definition}

\newtheorem*{déf*}{Définition} 

\newtheorem*{rem*}{Remarque} 
\newtheorem*{rems*}{Remarques} 

\newtheorem*{ex*}{Exemple} 
\newtheorem*{exs*}{Exemples} 

\newtheorem*{exo*}{Exercice} 
\newtheorem*{exos*}{Exercices} 

\theoremstyle{remark}

\newcommand{\disp}{\displaystyle}

\renewcommand{\a}{\alpha}

\renewcommand{\c}{\gamma}

\renewcommand{\d}{\delta}
\newcommand{\D}{\Delta}

\newcommand{\f}{\varphi}
\newcommand{\F}{\Phi}

\renewcommand{\k}{\kappa}
\renewcommand{\l}{\lambda}
\renewcommand{\L}{\Lambda}
\newcommand{\m}{\mu}

\newcommand{\om}{\omega}

\newcommand{\s}{\sigma}
\renewcommand{\S}{\Sigma}
\renewcommand{\t}{\theta}

\newcommand{\cC}{\mathcal{C}}

\newcommand{\cL}{\mathcal{L}}

\newcommand{\cP}{\mathcal{P}}
\newcommand{\cQ}{\mathcal{Q}}

\newcommand{\cS}{\mathcal{S}}
\newcommand{\cT}{\mathcal{T}}
\newcommand{\cU}{\mathcal{U}}
\newcommand{\cV}{\mathcal{V}}

\newcommand{\cZ}{\mathcal{Z}}

\newcommand{\as}{\ \mbox{\raisebox{.085ex}{$:$}\!$=$} \ } 

\newcommand{\st}{~|~} 

\newcommand{\End}[1]{\mathrm{End} \! \left( #1 \right)} 

\newcommand{\llist}[3]{#1_{#2} , \ldots , #1_{#3}} 

\newcommand{\inc}{\subset} 
\newcommand{\setmin}{\raisebox{0.45ex}{\scriptsize $\smallsetminus$}} 
\newcommand{\cart}{\! \times \!} 

\renewcommand{\to}{\longrightarrow} 
 
\newcommand{\I}[1]{\mathrm{I}_{\scriptscriptstyle #1}} 
\newcommand{\rest}[2]{#1_{\mathbf{|} #2}} 

\newcommand{\RR}{\mathbf{R}} 

\newcommand{\Rn}[1]{\mathbf{R}^{#1}} 
\newcommand{\Hn}[1]{\mathbf{H}^{#1}} 

\newcommand{\intr}[1]{\overset{\; _{\circ}}{#1}} 
\newcommand{\clos}[1]{\overline{#1}} 

\renewcommand{\leq}{\leqslant} 
\renewcommand{\geq}{\geqslant} 

\DeclareMathOperator{\atanh}{atanh}

\let\oldint\int
\renewcommand{\int}[4]{\oldint_{\! #1}^{#2} \!\!\! #3 \mathrm{d} #4} 
 
\renewcommand{\ll}[1]{\ell^{#1}} 

\renewcommand{\bar}{\overline} 

\let\olddet\det
\renewcommand{\det}[1]{\olddet{\! \left( #1 \right)}} 
\newcommand{\detb}[2]{\olddet\nolimits_{\! #1 \,}{\!\! #2}} 

\newcommand{\norm}[1]{\left\| #1 \right\|} 
\newcommand{\supnorm}[1]{\left\| #1 \right\|_{_{\infty}}} 
\newcommand{\opnorm}[1]{\left| \! \left| \! \left| #1 \right| \! \right| \! \right|} 

\newcommand{\diam}{\mathrm{diam}} 

\newcommand{\der}[2]{\frac{\mathrm{d} #1}{\mathrm{d} #2}} 
\newcommand{\Tg}[3]{T_{\! #1} #2 \! \cdot \! #3} 

\newcommand{\Cl}[1]{\mathrm{C}^{#1}} 

